# Introduction to Quantum Groups and Yang-Baxter Equation For Probabilists

Jeffrey Kuan

Accessibility statement: this PDF meets the technical standards of WCAG2.1AA, which complies with Ohio Administrative Policy IT-09 , Texas Administrative Code 206.70 and Title II of the Americans with Disabilities Act (effective April 24, 2026) . For the technical details related to screen reader usage, please visit my accessibility guidelines. For any issues, please send me an email.

A webpage version of this PDF, typeset in MathML, is also available. To block web crawlers, the webpage is password protected. The password is TaySwift13. As an additional benefit, the webpage will also have space for public comments and a list of updated errata, without the need to update the arXiv version every time an error is found.

## Abstract

Jeffrey Kuan

December 5, 2025

These are a set of lecture notes for a mini-course I gave at The University of Warwick from October 30th to November 1st, 2024. Recordings of the lectures are available on Oleg Zaboronski's webpage. The main body of the notes covers the content of the lectures, and provides an introduction to Drinfel'd-Jimbo quantum groups and the Yang-Baxter equation, with a probabilist as the target audience. The appendix contains several topics, requested by colleagues during my visit to the United Kingdom, which all depend on the main set of notes.

The notes begin by defining what it means for the asymmetric simple exclusion process (ASEP) to be integrable, in the sense of satisfying the Yang-Baxter equation. It then provides the algebraic background for the Yang-Baxter equation, by defining Drinfel'd-Jimbo groups as a quasi-triangular Hopf algebra. The algebraic background motivates generalizations of ASEP to stochastic vertex models and ``fused'' models. Each section corresponds to approximately an hour of lecture time.

The appendix covers the F.R.T. construction, Hecke algebras, the matrix product ansatz, and orthogonal polynomial vertex weights. The topics in the appendix can be read independently of each other.

Acknowledgements: Funding was provided by The London Mathematical Society and Lancaster University, University of Bristol, University of Edinburgh, Warwick University, Queen Mary University of London, and Ximera. These lecture notes were created with the help of Ximera, an interactive textbook platform hosted by The Ohio State University. The

Ximera Project is funded 2024-2026 (with no other external funding) by a $2,125,000 Open Textbooks Pilot Program grant from the federal Department of Education. I would like to thank Oleg Zaboronski for a careful reading of these notes.

# 1 Introduction

Broadly speaking, **integrable probability** is a branch of probability theory which studies models which have exact solutions. For this reason, these models are often called **exactly solvable.** From these exact solutions, one can derive precise **universal** asymptotics, such as the famed Tracy–Widom distribution. The origin of these solutions is usually due to algebraic symmetries underlying the model. In this set of lecture notes, we will introduce the relevant algebraic background with a probabilistic researcher as the target audience.

We first summarize some mathematical results since the early 1990s which contextualize the field in its current (2024) state. The well–known Tracy–Widom distribution was first discovered in [TW93], [TW94], where it describes the asymptotic fluctuations of the largest eigenvalues of the Gaussian Unitary Ensemble. In those papers, the distribution is defined in terms of a Fredholm determinant. Several years later, a surprising result of [BDJ99] proved that the Tracy–Widom distribution also occurs in the asymptotic fluctuations of longest increasing subsequences of random permutations. At the time, there was no explanation for why the same distribution would occur in the asymptotics of seemingly unrelated models.

The aforementioned exact solutions are a common theme in models with Tracy–Widom asymptotics. At this point in time, a complete list of all such models and relevant papers would be incredibly long, but for the purposes of these notes we can restrict to models with known quantum group symmetry. The canonical such model is the asymmetric simple exclusion process in one dimension. For an example of what these exact formulas look like, consider the result of [TW08], which shows that for an ASEP with initial conditions $y_1 < \cdots < y_N$, the probability at being in state $x_1 < \cdots < x_N$ at time $t$ is given by a sum over $N!$ terms, with each summand being a $N$–fold integral:

$$\sum_{\sigma \in S_N} \left(\frac{1}{2\pi i}\right)^N \int_{C_r} \cdots \int_{C_r} A_\sigma \prod_j \xi_{\sigma(j)}^{x_j - y_{\sigma(j)} - 1} e^{\sum_j \epsilon(\xi_j) t} d\xi_j.$$

For completeness, we state what these terms are. As usual $S_N$ is the symmetric group on $N$ elements, and an inversion of a permutation $\sigma$ is an ordered pair $(\sigma(j), \sigma(k))$ such that $j < k$ and $\sigma(j) > \sigma(k)$. The term $A_\sigma$ is defined to be

$$A_\sigma = \prod_{(\alpha,\beta) \text{ is an inversion of } \sigma} S_{\alpha\beta}$$

and $S_{\alpha\beta}$ is the "scattering term"

$$S_{\alpha\beta} = -\frac{1 + q\xi_\alpha \xi_\beta - (1+q)\xi_\alpha}{1 + q\xi_\alpha \xi_\beta - (1+q)\xi_\beta}.$$

The term $q$ is the asymmetry parameter. When $q = 0$, the totally asymmetric case, the formula reduces to a determinantal expression found in [Sch97a]. In this sense, the summation over the symmetric group $S_N$ is related to Leibniz's determinantal formula.

The asymptotic analysis proving Tracy–Widom convergence was done in the follow-up paper [TW09], generalizing the result for the totally asymmetric case ($q = 0$) in [Joh00]. Generally speaking, exact formulas and asymptotic analysis are easier when $q = 0$. However, we will later see that this is not true when using quantum groups.

At this point, the reader may be curious why any exact formulas should exist in the first place. This is due to the underlying quantum group symmetry of ASEP. To quote [TW08],

> "… it has been known for some time that the generator of ASEP is a similarity transformation of the quantum spin chain Hamiltonian known as the XXZ model. Since the XXZ Hamiltonian is diagonalizable by Bethe Ansatz, it is reasonable to expect that these ideas are useful for ASEP."

These notes will explicitly explain the underlying quantum group symmetry.

For even more historical context, the terminology for "integrable probability" originates in integrable systems from Hamiltonian mechanics. In Hamiltonian mechanics, the phase space is represented as a smooth manifold with even dimension $2n$, with coordinates denoted $q_1, \ldots, q_n, p_1, \ldots, p_n$ for the momentum and position. An integrable system is a system with $n$ conserved quantities which are mutually Poisson commuting, and the Liouville–Arnold theorem states that the equations of motion can be solved in quadratures. For an explicit example, the harmonic oscillator in one dimension (imagine an object attached to a frictionless spring) is integrable, and the conserved quantity is the total energy. (The appendix discusses this in more detail). In contrast, the three–body problem is not integrable, and its solutions are notoriously non–exact.

During the 1980s, many Soviet mathematicians and physicists introduced quantum mechanics into integrable systems. In quantum mechanics, quantities such as position, momentum and energy became operators which generally do not commute; additionally, the values of these quantities only take discrete, "quantized" values. In this context, the concept of "conserved quantities" becomes "commuting operators," and values of the quantities are eigenvalues of eigenstates. At the time, their approaches were called the "quantum inverse scattering method." Since then, the method has been generalized to abstract algebraic objects, such as Hopf algebras.

These notes will introduce one such algebraic object, known as **quantum groups,** with a particular focus on the Yang–Baxter equation. The exposition will use probability and mathematical physics as a motivation. The goal is for a reader with a probability background to be able to read contemporary (as of 2024) research papers in integrable probability.

# 2 ASEP, X.X.Z, and $U_q(sl_2)$

To begin to motivate the notes, we first introduce the asymmetric simple exclusion process (ASEP) and its relationship to the Heisenberg X.X.Z model and the quantum group $U_q(sl_2)$.

## 2.1 Definition of ASEP

ASEP was introduced in the mathematics community in [Spi70], and previously in a biology paper [MGP68]. In ASEP, particles randomly jump on a lattice, which we assume to be one–dimensional. At most one particle may occupy a site, and jumps to occupied sites are blocked (hence the term "exclusion"). Jumps are nearest neighbour (hence the term "simple"). If the jumps are continuous–time exponential clocks with left rates $\alpha$ and right rates $\beta$, then let $q = \sqrt{\beta/\alpha}$ denote the asymmetry parameter. If $q \neq 1$, then the model is asymmetric (sometimes called partially asymmetric), while for $q = 1$ the model is symmetric. For $q = 0$ or $q = \infty$ the model is called totally asymmetric. The symmetric exclusion process (without the word simple) can be defined more generally on an arbitrary graph [Spi70]. In principle, so can the asymmetric exclusion process, although this is not as well studied.

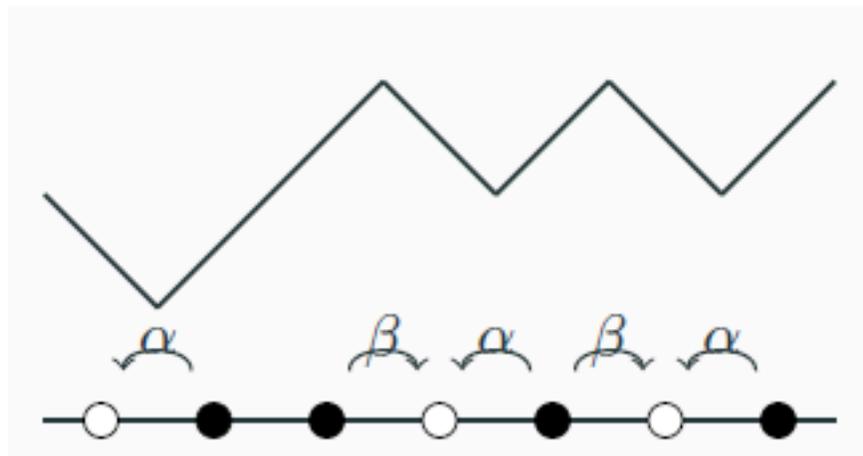

A particle configuration on a one-dimensional lattice on a line. Above each particle is a up-right path and above each non-particle is a down-right path. The paths connect to form a continuous jagged path. To the left of each particle is an arrow labeled α pointing to the lattice site to the left, and to the right of each particle is an arrow labeled β pointing to the lattice site to the right. However, if two particles are adjacent then there is no arrow between them.

Figure 1: An ASEP configuration with jump rates and the height function.

The generator of the simple exclusion process can be explicitly written. In the most elementary case where there are two lattice sites, then there are four possible configurations. Associate to each particle the vector $[1\ 0]$ and to each hole (i.e. a non–particle) the vector $[0\ 1]$. Tensoring the vectors together, we can associate to each configuration a canonical basis element of the four–dimensional vector space $\mathbb{C} \otimes \mathbb{C}$. The field here is chosen to be the complex numbers because it is algebraically closed, although of course probabilities are real numbers.

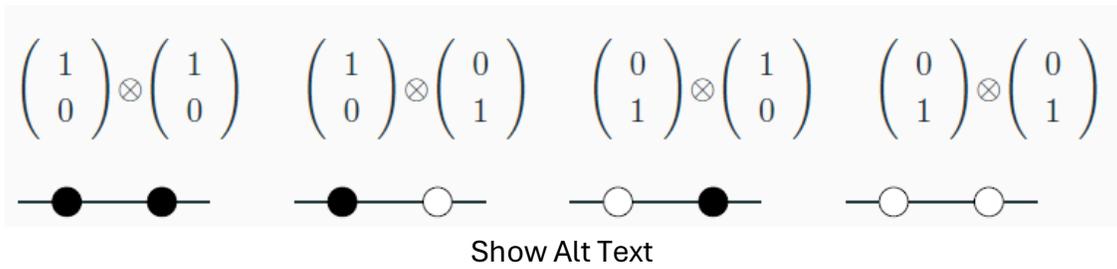

The four possible configurations of ASEP on two lattice sites. Above each particle is the vector $\begin{pmatrix}1\\0\end{pmatrix}$ and above each non-particle is the vector $\begin{pmatrix}0\\1\end{pmatrix}$. The two vectors have a tensor product sign between them.

Figure 2: Configurations of ASEP and their corresponding vectors.

With that set up, the generator is then a $4 \times 4$ matrix

$$\alpha \begin{pmatrix} 0 & 0 & 0 & 0 \\ 0 & -1 & 1 & 0 \\ 0 & q^2 & -q^2 & 0 \\ 0 & 0 & 0 & 0 \end{pmatrix}.$$

The constant $\alpha$ can be viewed as a time rescaling, so it can be removed without loss of generality.

In a more general setting, where there are $N$ lattice sites, the generator can be defined from the above $4 \times 4$ matrix, which we now denote $\mathcal{L}$. The generator is now a $2^N \otimes 2^N$ and acts on $(\mathbb{C}^2)^{\otimes N}$. Define $\mathcal{L}_{i,i+1}$ by

$$\mathcal{L}_{i,i+1} = (Id_2)^{i-1} \otimes \mathcal{L} \otimes (Id_2)^{N-1}.$$

With this notation, the generator is

$$\sum_{i=1}^{N-1} \mathcal{L}_{i,i+1}.$$

Similar notation using subscripts will be used throughout these notes.

## 2.2 ASEP and Yang–Baxter Equation

The corresponding notion of integrability in the quantum case comes from the Yang–Baxter equation. It was shown by Baxter [Bax72] that when a matrix solves the Yang–Baxter equation, there are corresponding "commuting transfer matrices," which are analogous to conserved quantities in the classical setting. In this section we will see that the generator

of ASEP satisfies the Yang–Baxter equation, justifying the description of ASEP as "integrable."

Using the same subscript notation as in the previous section, we have this definition:

*Definition*

We say that a matrix $R$ solves the braided Yang–Baxter Equation if:

$$R_{1,1}R_{2,3}R_{1,2} = R_{2,3}R_{1,2}R_{2,3}.$$

Somewhat confusingly, this equation is also sometimes called the braid equation or just the Yang–Baxter equation. An equivalent formulation is

$$\check{R}_{1,2}\check{R}_{1,3}\check{R}_{2,3} = \check{R}_{2,3}\check{R}_{1,3}\check{R}_{1,2}$$

where $\check{R} = P \circ R$ and $P(x \otimes y) = y \otimes x$ permutes the tensor power (exercise left to the reader).

A visual representation of Y.B.E is given by this image:

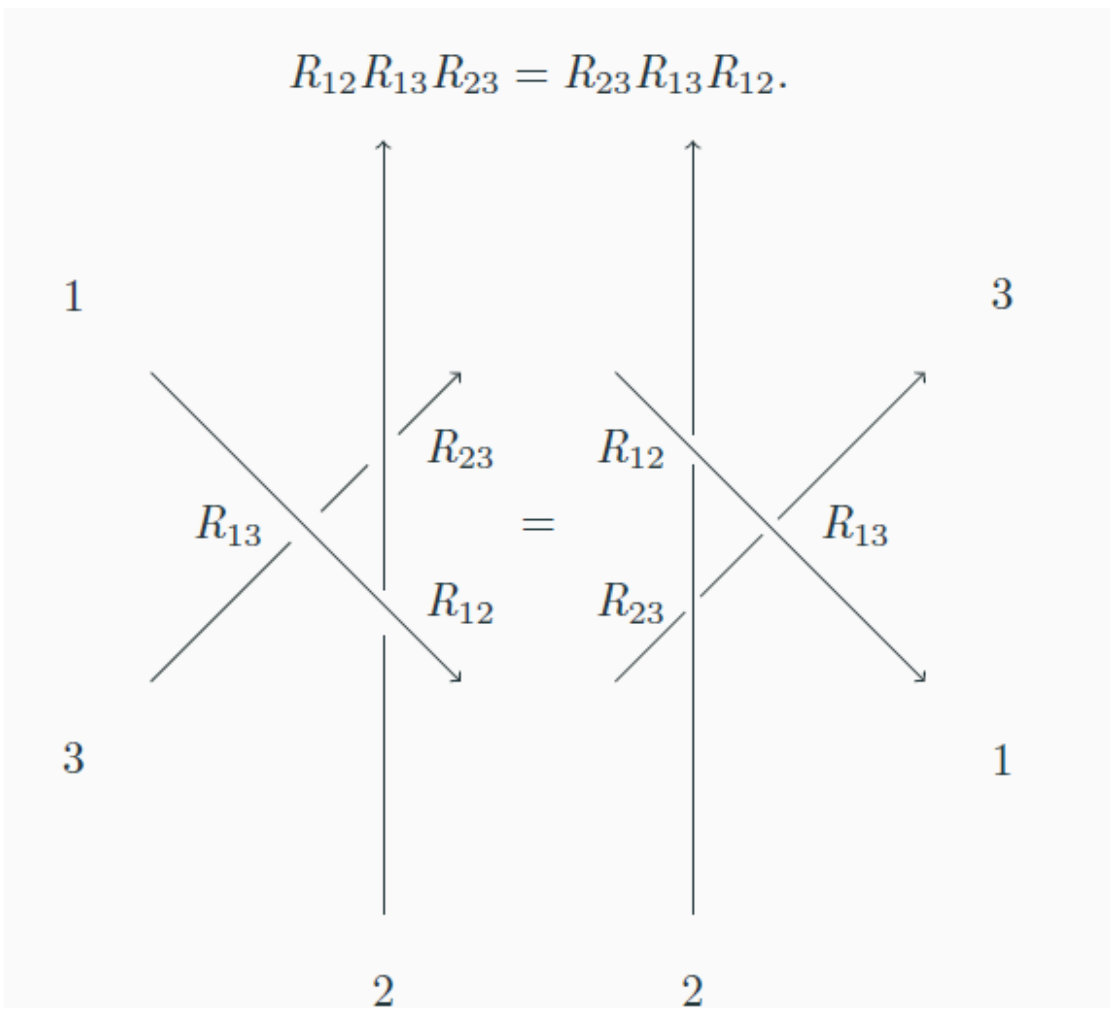

*Figure 3: The Yang-Baxter Equation.*

In some settings, the Y.B.E occurs in the scattering of particles with quantum mechanics considerations, such as the quantum inverse scattering method (Q.I.S.M.).

Given an arbitrary matrix $R$, one can check that it satisfies the Y.B.E through direct computation. However, it may be more helpful to consider some simple examples first.

If $P: V \otimes V \to V \otimes V$ denotes the permutation operator which sends $u \otimes v$ to $v \otimes u$, then the braided Yang–Baxter equation becomes an identity of transpositions:

$$(1\ 2) \circ (2\ 3) \circ (1\ 2) = (2\ 3) \circ (1\ 2) \circ (2\ 3).$$

This identity holds, since both sides equal the permutation $(1\ 3)$. An even more simple solution is the identity matrix.

A probabilist may remark that a permutation matrix and an identity matrix are both examples of a stochastic matrix, albeit somewhat trivial examples. It may then be natural to try to find more general stochastic matrices which solve the Yang–Baxter equation. The "simplest" generalization occurs when $V$ is two–dimensional, with basis $e_1, e_2$. Using exclusion processes as a prototype, we can define the operator $R_{\alpha\beta}$ by

$$R_{\alpha\beta}(e_1 \otimes e_1) = e_1 \otimes e_1, \quad R_{\alpha\beta}(e_2 \otimes e_2) = e_2 \otimes e_2$$

and

$$R_{\alpha\beta}(e_1 \otimes e_2) = (1 - \alpha)e_1 \otimes e_2 + \alpha e_2 \otimes e_1, R_{\alpha\beta}(e_2 \otimes e_1) = \beta e_1 \otimes e_2 + (1 - \beta)e_2 \otimes e_1.$$

Inserting this matrix into the Yang–Baxter equation and evaluating at $e_1 \otimes e_1 \otimes e_2$, one finds that $\alpha(1 - \beta) = 0$ is a necessary requirement for a solution to Y.B.E. Evaluating $R_{\alpha,1}$ and $R_{0,\beta}$ at $e_1 \otimes e_2 \otimes e_1$ and $e_2 \otimes e_1 \otimes e_1$ yields no additional equations. (Exercise left to the reader)

One could also evaluate Y.B.E at the three vectors $e_2 \otimes e_1 \otimes e_1, e_1 \otimes e_2 \otimes e_1$ and $e_1 \otimes e_1 \otimes e_2$, but instead we use another symmetry of ASEP. This is called particle–hole involution. In words, in an ASEP with drift to the right, the holes evolve as an ASEP with drift to the left. Symbolically, we define an involution $T$ of $V$ which switches $e_1$ and $e_2$. Then

$$T^{\otimes 2}R_{\alpha,\beta}T^{\otimes 2} = R_{\beta,\alpha}, \qquad T^{\otimes 2}R_{\alpha,\beta} = R_{1-\alpha,1-\beta}, \qquad R_{\alpha,\beta}T^{\otimes 2} = R_{1-\beta,1-\alpha}.$$

A mildly helpful observation here is that $T$ is stochastic, so its composition with other stochastic matrices is also stochastic. Using these identities, one immediately verifies Y.B.E for $R_{\alpha,0}$ and $R_{1,\beta}$. Furthermore, these two solutions are related up to the particle–hole involution.

One notices that $P \circ R_{\alpha,0} - id_4$ is the generator of ASEP at two lattice sites. It is then natural to ask whether or not addition by constants affects integrability. The answer is "no" (otherwise integrable probability would not exist). As some intuition, let us suppose that the generator has an explicit set of eigenvectors and eigenvalues. In general, an arbitrary Markov process will have complicated eigenvectors and eigenvalues, but for integrable models we expect "nice" eigenvectors and eigenvalues, from which probabilistic information can be extracted. Then addition by constants leaves the eigenvectors unaffected and only shifts the eigenvalues, leaving them just as "nice." Thus, it is fair to say that ASEP is integrable, at least on two lattice sites.

At this juncture, it is reasonable to object that nothing has been said about the integrability of ASEP on arbitrarily many lattice sites. We will put this question on hold, as the next section's discussion of the X.X.Z model will lay the foundation for answering this question.

## 2.3 X.X.Z. model

Having established that ASEP is an integrable model at two lattice sites, it is natural to ask about other integrable models that look similar to ASEP, so that we can avoid re–inventing the wheel. Thus we now turn our attention to a related model, called the quantum Heisenberg model. These were introduced by Werner Heisenberg [Hei28] to incorporate quantum mechanics into magnetism, and [GS92, ADHR94] pointed out the the X.X.Z Hamiltonian is conjugate to the ASEP generator.

Each lattice site has a microscopic magnetic dipole, which can either be up or down. If there are $N$ lattice sites, then the Hamiltonian acts on the $2^N$–dimensional vector space $(\mathbb{C}^2)^{\otimes N}$. So, at the very least, the state space of the X.X.Z model is the same as the ASEP state space.

To define the model, first recall the Pauli spin matrices:

*Definition*

The Pauli spin matrices are:

$$\sigma^1 = \begin{pmatrix} 0 & 1 \\ 1 & 0 \end{pmatrix}, \qquad \sigma^2 = \begin{pmatrix} 0 & -i \\ i & 0 \end{pmatrix}, \qquad \sigma^3 = \begin{pmatrix} 1 & 0 \\ 0 & -1 \end{pmatrix}.$$

The notation $\sigma^x, \sigma^y, \sigma^z$ are often used as well.

Note that these are traceless and are Hermitian, the latter property being important in physics (since Hermitian matrices have real eigenvalues, and these eigenvalues correspond to physical quantities). The Hamiltonian of the X.X.Z. model are defined in terms of the Pauli matrices:

*Definition*

Letting

$$\sigma_j^a = (Id)^{\otimes j-1} \otimes \sigma^a \otimes (Id)^{\otimes N-j}$$

so that $\sigma_j^a$ acts on the $j$–th lattice site while fixing the others, the Hamiltonian is then

$$\frac{-1}{2} \sum_{j=1}^{N} \left( J_x \sigma_j^1 \sigma_{j+1}^1 + J_y \sigma_j^2 \sigma_{j+1}^2 + J_z \sigma_j^3 \sigma_{j+1}^3 - h\sigma_j^3 \right).$$

where $J_x, J_y, J_z$ are some coupling constants and $h$ is the external field. We impose periodic boundary conditions, so that $\sigma_{N+1}^a = \sigma_1^a$. If $J_x, J_y, J_z$ are distinct then one obtains the X.Y.Z model. If $J_x = J_y \neq J_z$ then this is the X.X.Z model; if $J_x = J_y = J_z$ then it is the XXX model.

Each summand in the Hamiltonian can be written as a four by four matrix:

$$\begin{pmatrix} J_z + 2h & 0 & 0 & J_x - J_y \\ 0 & -J_z & J_x + J_y & 0 \\ 0 & J_x + J_y & -J_z & 0 \\ J_x - J_y & 0 & 0 & J_z - 2h \end{pmatrix}$$

If $J_x = J_y$, as assumed in the X.X.Z model, then the top–right and bottom–left entries are zero. In this case, the non–zero entries are at least the same as in the ASEP generator. Recalling the previous section's heuristic that adding by constants does not change integrability, we can set $J_z = 1, h = 0$ to obtain the matrix

$$Id_4 + \begin{pmatrix} 0 & 0 & 0 & 0 \\ 0 & -2 & 2J_x & 0 \\ 0 & 2J_x & -2 & 0 \\ 0 & 0 & 0 & 0 \end{pmatrix}.$$

At this point, it is worth noticing that the Y.B.E is a statement about operators on vector spaces, which should not depend on the choice of bases on the vector spaces. In other words, conjugation (which corresponds to a change of basis) should not affect integrability. As it turns out, conjugation by

$$\begin{pmatrix} 1 & 0 & 0 & 0 \\ 0 & 1 & \gamma - 1 & 0 \\ 0 & 0 & \gamma & 0 \\ 0 & 0 & 0 & 1 \end{pmatrix}$$

results in the ASEP generator on two sites. (Exercise left to the reader).

## 2.4 $sl_2$ symmetry

Next, we introduce the $sl_2$ symmetry of ASEP. This was first observed in the seminal paper of [Sch97].

So far, we have determined that ASEP is integrable and is conjugate to the X.X.Z Hamiltonian (assuming that we are ignoring the question of arbitrary lattice sizes). As it turns out, the Pauli matrices form a basis of the real Lie algebra $su_2$ corresponding to the real compact Lie group $SU(2)$. Here, a real Lie group means that the underlying field is the real numbers, not that all the entries are real (consider the unit circle in the complex plane as a real Lie group, for example).

The complexification of $SU(2)$ is the (non–compact) Lie group $SL(2, \mathbb{C})$ of traceless $2 \times 2$ matrices with complex entries. For a variety of reasons, some of which are historical, most mathematicians are introduced to the Lie algebra $sl_2$ rather than $su_2$, so we will use that notation.

Recall the following definition:

*Definition*

The Lie algebra $sl_2$ has basis denoted $e, f, h$ and brackets

$$[e, f] = h, \quad [h, e] = 2e, \quad [h, f] = -2f.$$

Note that the explicit matrices $E_{1,2}, E_{2,1}, E_{1,1} - E_{2,2}$ satisfy these relations, where $E_{ij}$ denotes the matrix with a 1 at the $i, j$–entry and a 0 everywhere else.

Returning to the XYZ model, one can calculate that the XXX model has $SU(2)$ symmetry (exercise to the reader) in the sense that

$$\sum_{j=1}^{N} \sigma_j^a$$

commutes with the Hamiltonian of the XXX model for $a = x, y, z$. However, for the X.X.Z model, the commutation only works for $a = z$. Therefore, whatever algebra is underlying the X.X.Z model (and hence that of ASEP), we know that it can not be as simple as $SU(2)$. Furthermore, since X.X.Z and ASEP have a free parameter, the algebra should be a one–parameter deformation of $SU(2)$. This leads into the next topic of quantum groups.

## 3 Quantum Group $U_q(sl_2)$

The Drinfel'd–Jimbo quantum groups [Dri85, Jim85] are quantizations of the Lie groups. Rather confusingly, the quantum groups are not themselves groups.

Before diving into the rigorous definition, one notes where the quantization occurs. Using probabilistic intuition, we expect that the SSEP should have an underlying $sl_2$ symmetry, whereas the ASEP should have the quantum group symmetry. However, SSEP and ASEP on a single lattice site are identical, so the quantization must require multiple lattice sites. We've established that in this case, the vector space is $(\mathbb{C}^{\otimes 2})^{\otimes N}$. This raises the question: how does $sl_2$ act on tensor powers?

To answer this question, we take an aside about Lie algebras, which are tangent spaces to Lie groups at the origin. What this implies is that for a Lie group $G$ with a Lie algebra $g$, we have that

$$X \in g \text{ implies } \exp(tX) \in G \text{ for all } t \in \mathbb{R}$$

and

$$X = \frac{d}{dt}\exp(tX)|_{t=0}.$$

Setting $g = g_X(t) = \exp(tX)$, the natural action of a group on tensor products is

$$g(v \otimes w) = gv \otimes gw.$$

Taking derivatives and heuristically using Leibniz's rule for derivative of products, we obtain

$$X(v \otimes w) = Xv \otimes w + v \otimes Xw.$$

This is expressed through the co–product $\Delta: g \to g \otimes g$ defined by

$$\Delta(X) = 1 \otimes X + X \otimes 1.$$

Since it is the co–product that allows for the algebra to act on multiple lattice sites (i.e. tensor powers), the quantization occurs in the co–product, not the algebra itself! As it turns out, the right deformation is

$$\Delta(e) = q^h \otimes e + e \otimes 1, \quad \Delta(f) = 1 \otimes f + f \otimes q^{-h}, \quad \Delta(h) = h \otimes 1 + 1 \otimes h, \quad (4)$$

where we interpret $q^h$ as a formal power series in $h$. This element is invertible with inverse denoted $q^{-h}$.

With these generators defined, the relations in the quantum group become  q h e q - h = q 2 e , q h f q - h = q - 2 f , e f - f e = q h - q - h q - q -1 (3)  Using the formal power series in h, these are actually the same relations as in the $q = 1$ case, where the last relation can be covered by formally using L'Hôpital's rule. More rigorously, there

Note that we do not deform the co–product of $h$, which is consistent with $\sigma^z$ commuting with the X.X.Z model.

Depending on the preferences of the author, the generator $h$ may be replaced with $k := q^h$ with the co–product becoming

$$\Delta(k) = k \otimes k.$$

(Exercise to the reader why that is the same co–product). A reader with algebraic background may recognize that $h$ is primitive and $k$ is group–like.

## 3.1 Hopf algebras

In this section, we will define Hopf algebras, which are named after Heinz Hopf. (This is a different Hopf than the namesake of the Cole–Hopf transformation, who is named Eberhard Hopf). The original paper possibly goes back to [8].

While I presented the co–product as the object that allows for multiple lattice sites, that is of course not how it is presented to algebraists. A more common description is that a co–algebra is dual to an algebra, with the co–product dual to the product. This subsection will provide a rigorous definition of a co–algebra.

One motivation is that the space of functions on an algebra is a co–algebra. The Lie algebra is actually itself a space of differential operators on a Lie group (it is the directional derivative), so is a space of functions on the algebra of smooth functions on the Lie group. Thus, you would expect a Lie algebra to have a co–algebra structure.

Recall that an algebra $A$ is a vector space with a multiplication operator

$$m: A \otimes A \to A$$

which is assumed to be bilinear. It is also associative in the sense that

$$m(x \otimes m(y \otimes z)) = m(m(x \otimes y) \otimes z),$$

which essentially means that $x(yz) = (xy)z$. This is equivalently stated as

$$m \circ (m \otimes id) = m \circ (id \otimes m).$$

Additionally, there is an identity element $1 \in A$ satisfying

$$m(1 \otimes x) = m(x \otimes 1) = x.$$

This can be written as an inclusion from the underlying field $F$ to the algebra $A$.

Before defining co–algebras, it will be useful to draw certain commutative diagrams. The associativity can be written as

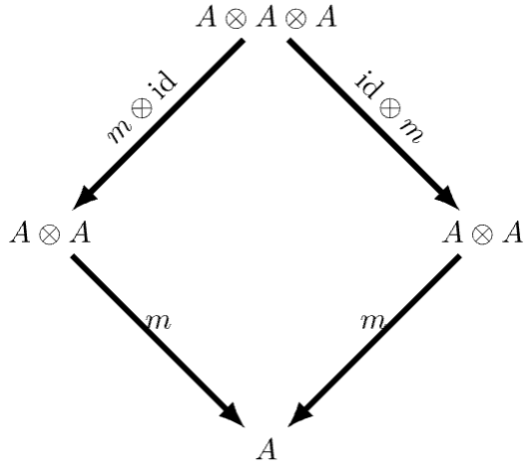

Show Alt Text

A commutative diagram, with four nodes and four arrows, in the shape of a diamond. The left side, from top to bottom, reads: $A \otimes A \otimes A \xrightarrow{m \otimes id} A \otimes A \xrightarrow{m} A$ and the right side, from top to bottom, reads: $A \otimes A \otimes A \xrightarrow{id \otimes m} A \otimes A \xrightarrow{m} A$

If one reverses the arrow of this diagram, which is also fortunately easy on Tikz, one obtains

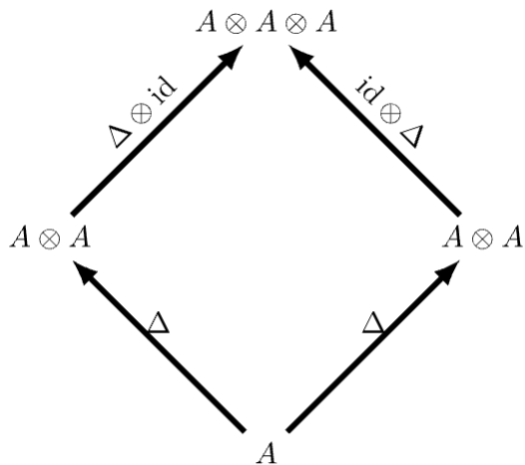

Show Alt Text

A commutative diagram, with four nodes and four arrows, in the shape of a diamond. The left side, from bottom to top, reads $A \xrightarrow{\Delta} A \otimes A \xrightarrow{\Delta \otimes id} A \otimes A \otimes A$ and the right side, from bottom to top, reads $A \xrightarrow{\Delta} A \otimes A \xrightarrow{id \otimes \Delta} A \otimes A \otimes A$

The condition of co–associativity means that this diagram commutes.

In a similar vein, the unit element can be understood as an inclusion from the base field to the algebra such that the following diagram commutes:

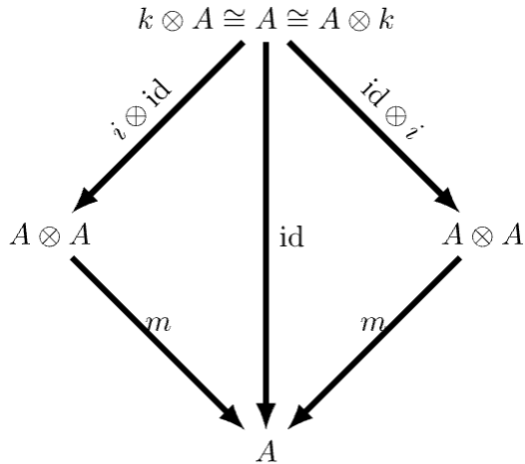

Show Alt Text

A commutative diagram, with four nodes and five arrows, in the shape of a diamond. From top to bottom, the left side reads $k \otimes A \cong A \cong A \otimes k \xrightarrow{\iota \otimes id} A \otimes A \xrightarrow{m} A$ and the right side reads $k \otimes A \cong A \cong A \otimes k \xrightarrow{id \otimes \iota} A \otimes A \xrightarrow{m} A$ and the middle goes from top to bottom and reads $k \otimes A \cong A \cong A \otimes k \xrightarrow{id} A$

Reversing the direction of the arrows, one obtains

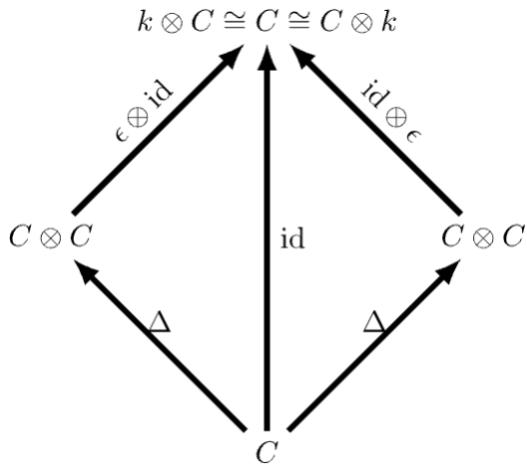

Show Alt Text

A commutative diagram, with four nodes and five arrows, in the shape of a diamond. From bottom to top, the left side reads $C \xrightarrow{\Delta} C \otimes C \xrightarrow{\epsilon \otimes id} k \otimes C \cong C \cong C \otimes k$ and the right side reads $C \xrightarrow{\Delta} C \otimes C \xrightarrow{id \otimes \epsilon} k \otimes C \cong C \cong C \otimes k$ and the middle goes from bottom to top and reads $C \xrightarrow{id} k \otimes C \cong C \cong C \otimes k$

A co–algebra is a vector space with a co–product and co-unit satisfying the above properties. For $U_q(\mathfrak{sl}_2)$ the co–unit is

$$\epsilon(e) = \epsilon(f) = 0, \quad \epsilon(q^h) = 1. \qquad (5)$$

A bi–algebra is an algebra and co–algebra such that the co–product and co–unit are also algebra homomorphisms.

A Hopf algebra is a bi–algebra with an anti–homomorphism $S$, called the underline{antipode}, such that this diagram commutes:

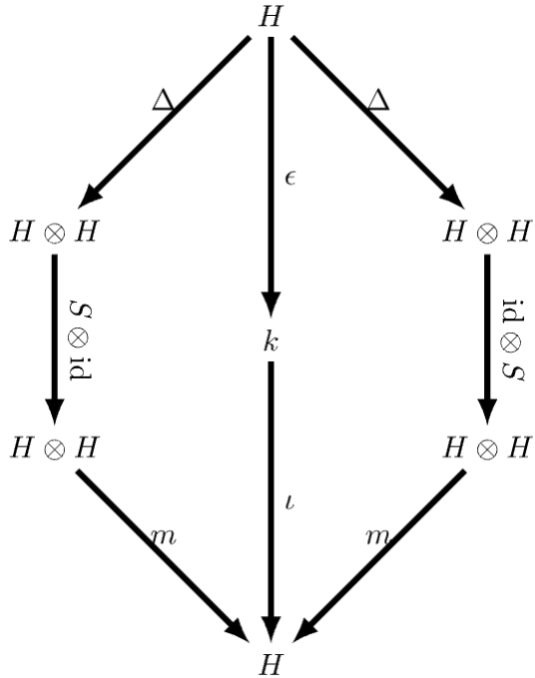

Show Alt Text

A commutative diagram with six nodes and eight arrows. It is in the shape of a hexagon, with three arrows on the left and three arrows on right. There are also two arrows going down the middle. From top to bottom, the left side reads $H \xrightarrow{\Delta} H \otimes H \xrightarrow{S \otimes id} H \otimes H \xrightarrow{m} H$ and the right side reads $H \xrightarrow{\Delta} H \otimes H \xrightarrow{id \otimes S} H \otimes H \xrightarrow{m} H$ and the middle goes from top to bottom and reads $H \xrightarrow{\epsilon} k \xrightarrow{\iota} H$

The antipode in this case is defined by

$$S(q^h) = q^{-h}$$
$$S(e) = -q^{-h}e$$
$$S(f) = -fq^h \qquad (6)$$

Finally, I get to $R$–matrices. A Hopf algebra $H$ is underline{co–commutative} if the coproduct $\Delta$ equals its reversal $\Delta =: = P \circ \Delta$ where $P: H \otimes H \to H \otimes H$ maps $x \otimes y$ to $y \otimes x$. A Hopf algebra is

<u>almost co–commutative</u> if there is an invertible element $R \in H \otimes H$ such that $R\Delta(x) = \Delta'(x)R$ for all $x \in H$.

An almost co–commutative Hopf algebra is <u>quasi–triangular</u> if the additional relations hold for the element $R$:

$$(\Delta \otimes id)(R) = R_{1,3}R_{2,3}, \quad (id \otimes \Delta)(R) = R_{1,3}R_{1,2}.$$

In a quasi–triangular Hopf algebra, the $R$–matrix satisfies the Yang–Baxter equation.

In [Dri87], the element R for $U_q(sl_2)$ is found:

$$q^{\frac{1}{2}h \otimes h} \sum_{i=0}^{\infty} (q - q^{-1})^i \frac{q^{-\frac{i(i+1)}{2}}}{[i]!} q^{\frac{i}{2}h} e^i \otimes q^{-\frac{i}{2}h} f^i. \quad (7)$$

Note that strictly speaking, since this is an infinite series, it is only an element in a certain completion. Indeed, Drinfeld's original formulation uses the algebra of formal power series $U[[\hbar]]$ where $\hbar$ denotes Planck's constant, and is related to q by $q = e^{\hbar}$. For probabilists, we can use usual analytic ideas of convergence, so this is less of an issue. The form of the element R above was expressed in [KS98Book]. The element R was recently used in [KZ24], so its expression is worth stating.

As a note on terminology: the algebraic element R is often called a <u>universal R-matrix</u> to contrast it with the usual R-matrix, which is actually a matrix. This leads to the somewhat humorous quip

> The definition of a universal R-matrix of a quantum group does not require either a matrix or a group.

We can thus summarize the definition of the quantum group as follows:

*Definition*

The Drinfel'd-Jimbo quantum group $U_q(sl_2)$ is the Hopf algebra with generators and relations defined by (3), with co-product defined by (4), with co-unit defined by (5), and with antipode defined by (6). With the element R defined by (7), it is a quasi-triangular Hopf algebra.

We can now summarize a main theorem.

*Theorem*

The asymmetric simple exclusion process is integrable, in the sense that its local generator on two lattice sites satisfies the Yang-Baxter equation, up to adding by constants times the identity matrix.

*Proof*

The quantum group $U_q(sl_2)$ is quasi-triangular as a Hopf algebra, so its R-matrix satisfies the Yang-Baxter equation. By taking its representation on a two-dimensional vector space, one obtains the ASEP generator, up to adding $cI$

## 4 Explicit solutions to Y.B.E

### 4.1 Affine Lie algebras and vertex models

Having shown that solutions to Yang–Baxter Equation can come from Hopf algebras, we now introduce a spectral–dependent Yang–Baxter Equation that comes from affine Lie algebras. This Lie algebra is now infinite–dimensional, which makes its representation theory trickier. The definition of affine Lie algebras is given below:

*Definition*

Given a finite–dimensional Lie algebra $\mathfrak{g}$, define its affinization by

$$\hat{\mathfrak{g}} = \mathfrak{g} \otimes_{\mathbb{C}} \oplus \mathbb{C}[t, t^{-1}] \oplus \mathbb{C}c$$

where $c$ is a central element. In words, this is a one–dimensional central extension of the loop algebra over $\mathfrak{g}$. The Lie bracket is then defined by

$$[a \otimes t^n, b \otimes t^m] = [a, b] \otimes t^{n+m} + (a, b)n\delta_{m,-n}c$$

where $(\cdot,\cdot)$ is the Cartan–Killing form on $\mathfrak{g}$.

There is an analogous quantum group $\mathcal{U}_q(\mathfrak{g})$, which has an $R$–matrix satisfying Yang–Baxter equation. For modules where $c$ acts as 0 (level 0 representations), we can define an <u>evaluation module</u> by taking a module of $\mathfrak{g}$ and then letting $t$ act as multiplication by $z$. If the element $R$ now acts on evaluation modules, then it depends on a parameter which we call the spectral parameter. This leads to the definition:

*Definition*

The spectral-dependent Yang-Baxter Equation is:

$$\check{R}_{1,2}(z)\check{R}_{1,3}(zw)\check{R}_{2,3}(w) = \check{R}_{2,3}(w)\check{R}_{1,3}(zw)\check{R}_{1,2}(z).$$

If we attempt to find $4 \times 4$ solutions to the spectral–dependent Y.B.E analogous to ASEP, we get a two–parameter family of solutions, which end up being the weights of the stochastic six vertex model, which we now define.

The stochastic six vertex model was defined in [GS92], with the generic six vertex model going back to [Pau35]. The six weights are displayed below:

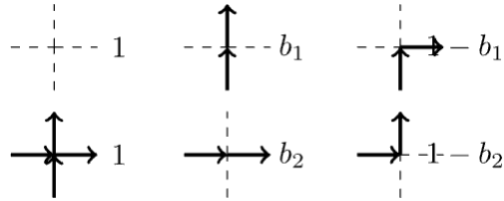

Show Alt Text

The six vertices of the six-vertex model. Each vertex is shaped like a plus sign with a certain number of arrows coming in from the bottom and left and the same number of arrows going out to the top and right. If there are no arrows coming in and no arrows going out, the weight is 1. If there are two arrows coming in and two arrows going out, the weight is 1. If exactly one arrow comes in from the bottom and goes up, the weight is $b_1$. If exactly one arrow comes in from the left and goes right, the weight is $b_2$. If exactly one arrow comes in from the bottom and goes right, the weight is $1 - b_1$. If exactly one arrows come in from the left and goes up, the weight is $1 - b_2$.

*Figure 4: The six weights of the stochastic six--vertex model.*

The number of arrows coming in from the bottom and left must equal the number of arrows coming out from the top and right. Thus, there are $6 = (4//2)$ possible vertices. On a two–dimensional lattice, perhaps with fixed boundary conditions, this defines a weight on configurations:

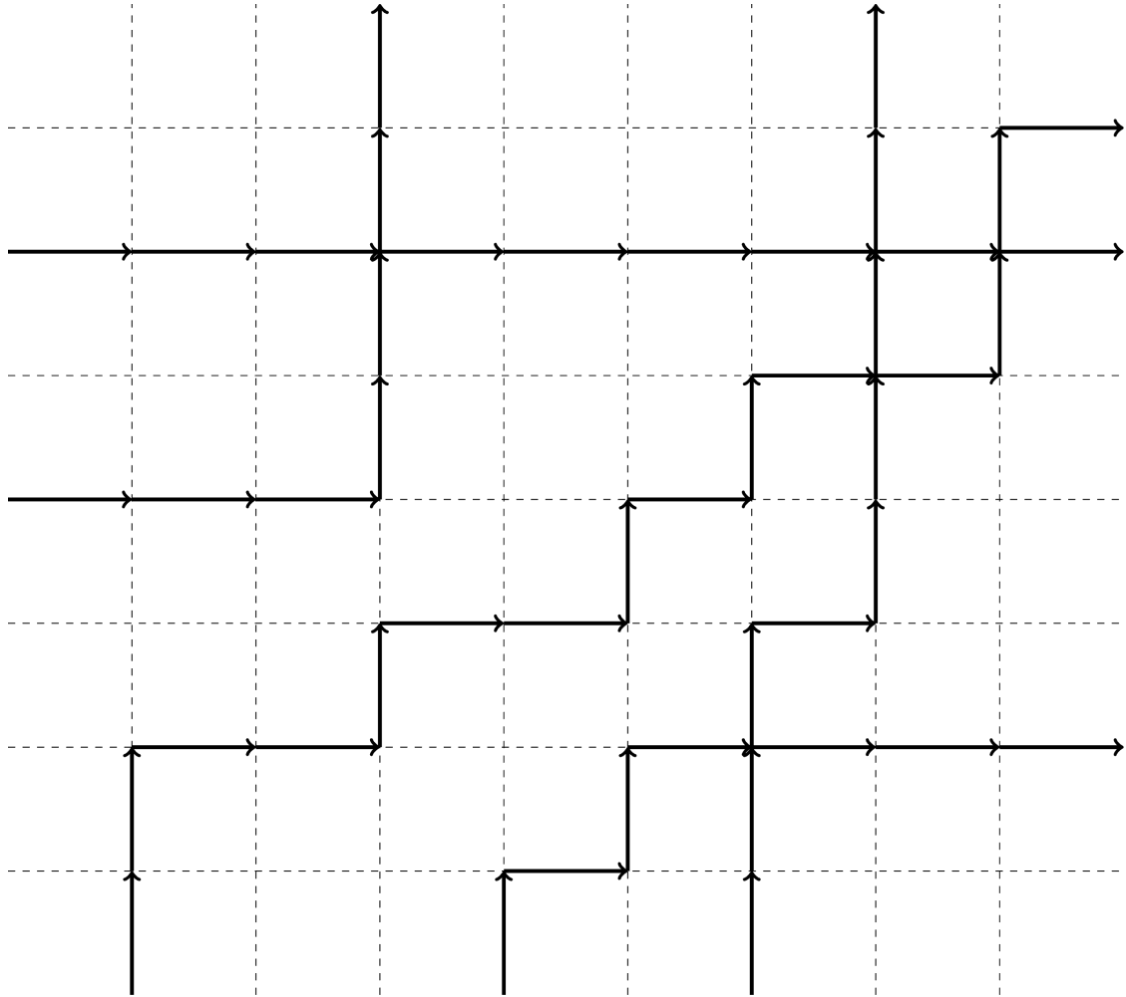

A large two-dimensional lattice where each lattice site has a vertex from the six-vertex model. It appears that the arrows form paths that zigzag up and to the right.

Figure 5: A sample of the stochastic six--vertex model on a two-dimensional lattice.

## 4.2 Fusion

In general, we do not really want to brute force solutions to the Yang–Baxter equation in arbitrary representations. Here, I present a way to construct new solutions from old solutions. The original idea comes from [CP16].

In linear algebra, for any vector space $V$, a symmetric tensor is an element of $V^{\otimes N}$ preserved by any $P_{i,i+1}$. For quantum groups, however, we have a $q$–deformation in the tensor powers so that $P$ and $\Delta$ no longer commute. However, there is still a symmetric tensor power

$$Sym^N V \subseteq V^{\otimes N}.$$

In lieu of an extensive algebraic discussion, we present it in probabilistic language. Recall Pitman–Rogers intertwining [RP81] (also known as lumping [KS60]): suppose there are two Markov kernels $\Lambda$ and $\Phi$ such that $\Lambda\Phi = Id$. Let $Q = \Lambda P \Phi$ and assume $\Lambda P = Q\Lambda$. Then $\Phi$ defines a Markov projection from $P$ to $Q$:

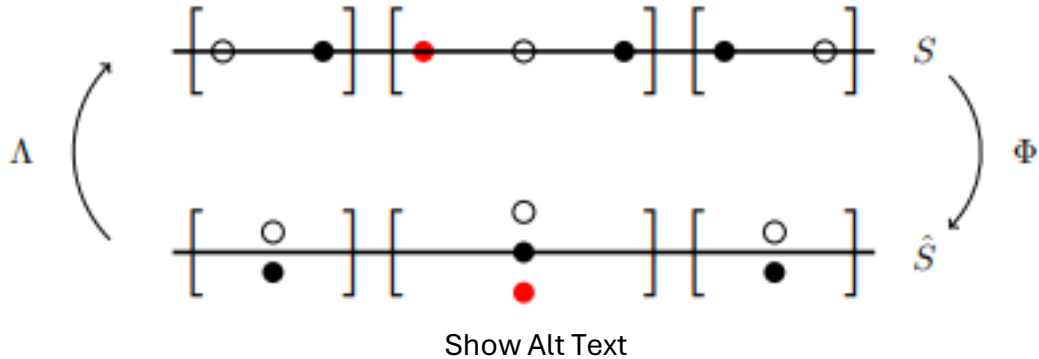

The top portion of the image shows 7 lattice sites with some particles, grouped into blocks of size 2,3 and 2. The bottom portion of the image shows three blocks, which can contain up to 2,3, or 2 particles. The particles in the top portion are mapped into the blocks on the bottom portion. The top portion is labeled S and the bottom portion is labeled $\hat{S}$. An arrow labeled $\Lambda$ goes from bottom to top and an arrow labeled $\Phi$ goes from top to bottom.

Figure 6: *A visualization of fusion of lattice sites.*

More specifically, if $X(t)$ is a Markov process on state space $S$ and $\Phi: S \to \hat{S}$ is deterministic, with a Markov kernel $\Lambda: \hat{S} \to S$, then $\Phi(X(t))$ is a Markov process on $\hat{S}$.

This procedure is sometimes called fusion.

This intertwining can be constructed algebraically, but instead we draw a metaphor with stationary measures. We want to show that

$$\Lambda P = \Lambda P \Phi \Lambda.$$

If for all $v$, we have that $v\Lambda$ is stationary under $P$, i.e. $v\Lambda P = v\Lambda$, then the above intertwining becomes

$$v\Lambda = v\Lambda\Phi\Lambda$$

which then follows from $\Lambda\Phi = Id$.

Calculating the fused vertex weights then requires orthogonal polynomials (please reference the theorem in Section 5.4).

At this point, a natural question arises: do these weights satisfy Y.B.E? For the stochastic vertex models, they do. For interacting particle systems such as ASEP, the situation seems to be dependent on the choice of Lie algebra.

## 4.3 Construction with Central Elements

Here, we mention a method for constructing non–integrable interacting particle systems using central elements. The original argument goes back to [Sch97] for ASEP (which happens to be integrable), and was then generalized in [CGRS16] for a partial exclusion process, where an arbitrary number of particles may occupy a site.

From here, let $k = q^h$. We define some subalgebras by letting $U^{>0}$ be generated by $e$, and $U^{<0}$ generated by $f$ and $U^0$ generated by $k$. Analogously, $U^{\geq 0}$ is generated by $e, k$ and $U^{\leq 0}$ is generated by $f, k$. There are central elements, sometimes called Casimir elements, that have the form

$$C \in U^0 + U^{<0}U^{>0}.$$

One reason for introducing these subalgebras is because of the highest weight representations. We call $V$ a highest weight representation if there exists a vector $v$, called the highest weight vector, such that

$$U^{>0}v = 0$$

and $V$ is generated by $U^{\leq 0}v$. Then the action of the Casimir on the highest weight vector is given by

$$Cv = av$$

for some complex number $a$. Also let $\{v_\mu : \mu \in \mathcal{B}\}$ denote a basis of $V$, where each $v_\mu$ can be expressed as $F^\mu v$ for some $F^\mu \in U^{\leq 0}$.

We define the Hamiltonian $H$ by taking the co–product of $C$ on tensor powers of representations:

$$H = \Delta(C)|_{V \otimes V}.$$

We want to find a diagonal matrix $G$ such that $G^{-1}HG - c\mathrm{Id}$ is the generator of a continuous–time Markov process. As a necessary condition, we need that $G^{-1}HG$ has rows which sum to $c$; non–negativity is more difficult to ensure. One way to achieve this is with an eigenvector of eigenvalue $c$. Namely, assume that $g$ is a vector such that $Hg = cg$. If $g$ has the form

$$g = \sum_{\boldsymbol{\mu} = (\mu_1, \mu_2) \in \mathcal{B} \times \mathcal{B}} g(\boldsymbol{\mu}) v_{\mu_1} \otimes v_{\mu_2}$$

then the diagonal matrix $G(\boldsymbol{\mu}, \boldsymbol{\mu}) = g(\boldsymbol{\mu})$ will satisfy the necessary requirement, as long as $g(\boldsymbol{\mu})$ is always nonzero.

As it turns out, $g = v \otimes v$ (the tensor power of the highest weight vector) will satisfy

$$Hg = a^2 g$$

due to the form of the co–product and the central element. However, this choice will have zero entries. So instead we take $g$ of the form

$$\sum_{\mu=(\mu_1,\mu_2)} c(\mu_1,\mu_2) \Delta(F^{\mu_1} F^{\mu_2}) v \otimes v$$

for some appropriate choice of coefficients $c(\mu_1, \mu_2)$ (reference the paper [Kua18] for more detailed calculations).

# 5 Appendix

## 5.1 F.R.T Construction

### 5.1.1 Motivation

Given a solution $R(\lambda)$ of the Yang–Baxter equation, the F.R.T construction of [FRT88] produces a quasi–triangular Hopf algebra. To begin, suppose suppose the $R$–matrix $R(\lambda)$ acts on the tensor power of vector spaces $V \otimes V$. Let $L(\lambda)$ be an operator acting on the tensor product $V \otimes W$. Now consider the equation

$$R_{1,2}(\lambda - \mu) L_{1,3}(\lambda) L_{2,3}(\mu) = L_{2,3}(\mu) L_{1,3}(\lambda) R_{1,2}(\lambda - \mu)$$

where both sides act on the tensor product $V \otimes V \otimes W$. One solution occurs when the vector spaces $V$ and $W$ are the same and $L$ is just equal to $R$. For the remainder of these notes, we consider the simpler case when $V = W$.

Before launching into the most generalize case, it will be enlightening to consider the case of $U_q(sl_2)$. Most of this portion of the exposition comes from the book [KS98Book]. Consider the $4 \times 4$ matrix

$$R := \begin{pmatrix} q^{-2} & 0 & 0 & 0 \\ 0 & 0 & q^{-1} & 0 \\ 0 & q^{-1} & q^{-2} - 1 & 0 \\ 0 & 0 & 0 & q^{-2} \end{pmatrix}.$$

We view this naturally as an operator on the tensor power $\mathbb{C}^2 \otimes \mathbb{C}^2$. One can verify that this matrix solves the braided Yang–Baxter equation

$$R_{1,2} R_{2,3} R_{1,2} = R_{2,3} R_{1,2} R_{2,3}$$

and

$$(R - q^{-2})(R + 1) = 0.$$

To relate this to $U_q(sl_2)$, we first define a representation $\rho_m$ on a $(m+1)$–dimensional vector space, with basis denoted $\{v_0, \ldots, v_m\}$. The representation is defined by (fixing a typo in [KS98Book])

$$\rho_m(e)v_k = \frac{q^{m-k} - q^{k-m}}{q - q^{-1}} v_{k+1},$$

$$\rho_m(f)v_k = \frac{q^k - q^{-k}}{q - q^{-1}} v_{k-1},$$

$$\rho_m(h)v_k = (2k - m)v_k.$$

For $i, j = 1, 2$ let $t_{i,j}$ denote the matrix entries of the two–dimensional representation $\rho_1$, and let $T$ denote the $2 \times 2$ matrix with entries $t_{i,j}$. In other words, we consider $t_{i,j}$ as abstract, non–commuting generators of the algebra $Hom(U_q(sl_2), \mathbb{C})$, defined by

$$t_{i,j}(u) = [\rho_1(u)]_{i,j}.$$

The product in dual bi–algebras is defined by the co–product in the original algebra. In other words,

$$(t_{i,j} t_{k,l})(u) = (t_{i,j} \otimes t_{k,l})(\Delta(u)).$$

Then one can verify that (fixing a typo in [KS98])

$$R(1 \otimes T)(T \otimes 1) = (T \otimes 1)(1 \otimes T)R.$$

The relations can be explicitly written (details omitted) as

$$t_{1,1}t_{1,2} = q^{-1}t_{1,2}t_{1,1}, \qquad t_{1,1}t_{2,1} = q^{-1}t_{2,1}t_{1,1},$$
$$t_{1,2}t_{2,2} = q^{-1}t_{2,2}t_{1,2}, \qquad t_{2,1}t_{2,2} = q^{-1}t_{2,2}t_{2,1},$$
$$t_{1,2}t_{2,1} = t_{2,1}t_{1,2}, \qquad t_{1,1}t_{2,2} - t_{2,2}t_{1,1} = (q^{-1} - q)t_{1,2}t_{2,1},$$
$$t_{1,1}t_{2,2} - q^{-1}t_{1,2}t_{2,1} = 1$$

Of course, writing in the RTT form is much more concise. If the last relation is removed, then we obtain the dual to $U_q(gl_2)$.

When $q = 1$, this reduces to the commutative algebra $C[SL(2)]$, the usual algebra of functions on the group $SL(2)$. This should not be surprising: the quantum group $U_q(sl_2)$ is non–co–commutative except when $q = 1$, so its dual should be non–commutative except when $q = 1$. The algebra defined above is denoted $C_q[SL(2)]$.

The co–product in a dual is defined from the product in the original algebra. In this case, this means that

$$\Delta(t_{i,j})(u_1 \otimes u_2) = t_{i,j}(u_1 u_2).$$

The explicit formula is then

$$\Delta(t_{i,j}) = \sum_k t_{i,k} \otimes t_{k,j}.$$

Note that there is no $q$ in this formula, which should not surprise us. Since the quantization in $U_q(sl_2)$ deforms the co–product and not the product, the corresponding deformation in the dual should not deform the co–product, but rather the product. In the more general case, we have this theorem:

### Theorem

Suppose that a matrix $R$ satisfies the braided Yang–Baxter equation,

$$R_{1,2}R_{2,3}R_{1,2} = R_{2,3}R_{1,2}R_{2,3}$$

then $T(R)$ is a bialgebra with generators $\{t_{i,j}\}$ and relations given by

$$R(1 \otimes T)(T \otimes 1) = (T \otimes 1)(1 \otimes T)R$$

and co–product

$$\Delta(t_{i,j}) = \sum_k t_{i,k} \otimes t_{k,j}.$$

To recover a quantum group analogous to a Drinfel'd–Jimbo quantum group, we would need to take the restricted dual, which is contained in $Hom(T, \mathbb{C})$. This is generated by all matrix entries of finite–dimensional modules. Since $T$ is generally infinite–dimensional, the full dual is too large to reasonably consider. We elaborate on this in the next section.

### 5.1.2 Dual algebras

Recall that the dual of a co–algebra is an algebra, and vice versa. The multiplication $m$ in the algebra is defined by the co–product $\Delta$ in the co–algebra. To describe this in more detail, we recall Sweedler's notation. For any $c$ in any co–algebra, we write

$$\Delta(c) = \sum_{(c)} c_{(1)} \otimes c_{(2)}.$$

Sometimes the summation is dropped, and in that case it is called sumless Sweedler's notation.

Now, let $A$ be an algebra over a field $k$ and let $C = Hom_k(A, k)$. Then the co–product is defined by

$$\Delta(c)(a_1 \otimes a_2) = \sum_{(c)} c_1(a_1)c_2(a_2) = c(a_1 a_2)$$

In our situation, we are considering Hopf algebras, which are both algebras and co–algebras. So in this, it is more convenient to use the notation

$$\langle c, a \rangle = c(a).$$

It would perhaps be more helpful to illustrate this with the quantum group $U_q(sl_2)$. More particularly, we have that for $u \in U_q(sl_2)$,

$$\langle t_{ij}, u \rangle = [\rho_1(u)]_{ij}.$$

In words, we say that $t_{ij}$ is the $(i,j)$ matrix entry in the fundamental representation. Suppose that we did not know the co–product applied to $e$. We do know that since $e$ increases weights by one (from classical Lie theory) that

$$\Delta(e) = k^{(1)} \otimes e + e \otimes k^{(2)}$$

for some $k^{(1)}, k^{(2)} \in U^0$. Then

$$\langle k^{(1)}, t_{i,j}\rangle\langle k^{(2)}, t_{i',j'}\rangle = \langle \Delta(e), t_{i,j} \otimes t_{i',j'}\rangle = \langle e, t_{i,j}t_{i',j'}\rangle.$$

We then need to solve for the eight variables $\{\langle k^{(a)}, t_{i,j}\rangle : a, i, j = 1,2\}$. The relations among $t_{i,j}$ give seven equations, and there is an additional equation because we know the co–product when $q = 1$. More explicitly, copying and pasting the relations, we have

$$\langle k^{(1)}, t_{1,1}\rangle\langle e, t_{1,2}\rangle + \langle e, t_{1,1}\rangle\langle k^{(2)}, t_{1,2}\rangle = q^{-1}\big(\langle k^{(1)}, t_{1,2}\rangle\langle e, t_{1,1}\rangle + \langle e, t_{1,2}\rangle\langle k^{(2)}, t_{1,1}\rangle\big)$$

$$\langle k^{(1)}, t_{1,1}\rangle\langle e, t_{2,1}\rangle + \langle e, t_{1,1}\rangle\langle k^{(2)}, t_{2,1}\rangle = q^{-1}\big(\langle k^{(1)}, t_{2,1}\rangle\langle e, t_{1,1}\rangle + \langle e, t_{2,1}\rangle\langle k^{(2)}, t_{1,1}\rangle\big)$$

$$\langle k^{(1)}, t_{1,2}\rangle\langle e, t_{2,2}\rangle + \langle e, t_{1,2}\rangle\langle k^{(2)}, t_{2,2}\rangle = q^{-1}\big(\langle k^{(1)}, t_{2,2}\rangle\langle e, t_{1,2}\rangle + \langle e, t_{2,2}\rangle\langle k^{(2)}, t_{1,2}\rangle\big)$$

$$\langle k^{(1)}, t_{2,1}\rangle\langle e, t_{2,2}\rangle + \langle e, t_{2,1}\rangle\langle k^{(2)}, t_{2,2}\rangle = q^{-1}\big(\langle k^{(1)}, t_{2,2}\rangle\langle e, t_{2,1}\rangle + \langle e, t_{2,2}\rangle\langle k^{(2)}, t_{2,1}\rangle\big)$$

$$\langle k^{(1)}, t_{1,2}\rangle\langle e, t_{2,1}\rangle + \langle e, t_{1,2}\rangle\langle k^{(2)}, t_{2,1}\rangle = \langle k^{(1)}, t_{2,1}\rangle\langle e, t_{1,2}\rangle + \langle e, t_{2,1}\rangle\langle k^{(2)}, t_{1,2}\rangle$$

$$\big(\langle k^{(1)}, t_{1,1}\rangle\langle e, t_{2,2}\rangle + \langle e, t_{1,1}\rangle\langle k^{(2)}, t_{2,2}\rangle\big) - \big(\langle k^{(1)}, t_{2,2}\rangle\langle e, t_{1,1}\rangle + \langle e, t_{2,2}\rangle\langle k^{(2)}, t_{1,1}\rangle\big)$$
$$= (q^{-1} - q)\big(\langle k^{(1)}, t_{1,2}\rangle\langle e, t_{2,1}\rangle + \langle e, t_{1,2}\rangle\langle k^{(2)}, t_{2,1}\rangle\big)$$

$$\big(\langle k^{(1)}, t_{1,1}\rangle\langle e, t_{2,2}\rangle + \langle e, t_{1,1}\rangle\langle k^{(2)}, t_{2,2}\rangle\big) - q^{-1}\big(\langle e, t_{1,2}\rangle\langle k^{(2)}, t_{2,1}\rangle + \langle k^{(1)}, t_{1,2}\rangle\langle e, t_{2,1}\rangle\big)$$
$$= \langle k^{(1)}, 1\rangle\langle e, 1\rangle + \langle e, 1\rangle\langle k^{(2)}, 1\rangle$$

Using that $\langle e, t_{1,2}\rangle = 1$ and $\langle e, t_{1,1}\rangle = \langle e, t_{2,1}\rangle = \langle e, t_{2,2}\rangle = 0$ we get that

$$\langle k^{(1)}, t_{1,1}\rangle = q^{-1}\langle k^{(2)}, t_{1,1}\rangle$$

$$\langle k^{(2)}, t_{2,2}\rangle = q^{-1}\langle k^{(1)}, t_{2,2}\rangle$$

This does not have a unique solution, but once $k^{(2)}$ is chosen then $k^{(1)}$ is determined. One solution is

$$k^{(2)} = 1, \quad k^{(1)} = k^{-1}.$$

The existence of more than one solution is why different papers have different conventions for the co–product.

### 5.2 Hecke Algebras

#### 5.2.1 Group Algebras and Coxeter Groups

For a finite group $G$ and a field $k$, the group algebra $k[G]$ is the algebra over $k$ with dimension $|G|$ and a basis indexed by the elements of $G$. The multiplication in $k[G]$ is given by the multiplication in $G$. As an example, the group algebra $\mathbb{C}[S_n]$ has dimension $n!$ with basis denoted by permutations $\sigma \in S_n$. The multiplication of two basis vectors equals another basis vector, induced from the multiplication in the symmetric group $S_n$.

The symmetric group $S_n$ is an example of a Coxeter group. We can present $S_n$ with generators $s_1, \ldots, s_{n-1}$ where each $s_i$ is the permutation $(i\, i+1)$. Then the relations are

$$s_i^2 = 1, \quad (s_i s_{i+1})^3 = 1, \quad (s_i s_j)^2 = 1 \text{ if } |i - j| \neq 1.$$

The condition $(s_i s_j)^2 = 1$ is equivalent to $s_i s_j = s_j^{-1} s_i^{-1} = s_j s_i$, which means that $s_i$ and $s_j$ commute. This generalizes in the following way:

### Definition

A **Coxeter group** is a group $W$ generated by a finite set $S = \{s_1, \ldots, s_n\}$ of involutions, subject only to relations of the form $(s_i s_j)^{m_{ij}} = 1$, where $m_{ii} = 1$ for all $i$, $m_{ij} = m_{ji} \geq 2$ for $i \neq j$, and $m_{ij} \in \mathbb{N} \cup \{\infty\}$ (with the convention that no relation is imposed when $m_{ij} = \infty$). The integers $m_{ij}$ can be collected in a **Coxeter matrix** $M = (m_{ij})_{1 \leq i,j \leq n}$.

For example, the symmetric group $S_{n+1}$ has Coxeter matrix

$$\begin{cases} m_{ii} = 1 \\ m_{i,i+1} = 2 \\ m_{ij} = 3, |i - j| > 1. \end{cases}$$

Of course, there are strict restrictions on which matrices can be Coxeter matrices. The classification of Coxeter matrices can be found in [Cox35], and most probabilists could use this classification as a black box.

### 5.2.2 Hecke Algebras

The Hecke algebra is a q–deformation of (the group algebra) of a Coxeter group, analogously to the quantum group being a q–deformation of (the universal enveloping algebra) of a Lie algebra. First, note that if $m_{ij} = 3$ then

$$s_i s_j s_i s_j s_i s_j = 1$$

can be rewritten as

$$s_i s_j s_i = s_j s_i s_j$$

since $s_i^2 = 1$ for all i. By rewriting this way, the equation can be justifiably called a braid relation.

### Definition

Given a Coxeter group $W$ with generating set $S$, the Hecke algebra $H_q(W)$ (or $H(W, q)$) is a $q$-deformation of the group algebra $\mathbb{C}[W]$. It is generated by elements $\{T_i : 1 \leq i \leq n\}$ satisfying:

- The braid relations: $(T_i T_j T_i \cdots) = (T_j T_i T_j \cdots)$ with $m_{ij}$ factors on each side, corresponding to the Coxeter relations, and
- The quadratic relation: $(T_i - q)(T_i + 1) = 0$ for each i.

When $q = 1$, the quadratic relation becomes $T_i^2 = 1$, and the Hecke algebra specializes to the group algebra $\mathbb{C}[W]$.

As an example of a Hecke algebra, we have that $H_q(S_{n+1})$ is generated by $T_1, \ldots, T_n$ with relations

$$(T_i - q)(T_i + 1) = 0, \quad T_i T_{i+1} T_i = T_{i+1} T_i T_{i+1}.$$

For the type B Hecke algebra, which deforms the hyperoctohedral group, there is in fact a two–parameter deformation. This algebra is generated by $T_0, T_1, \ldots, T_{n-1}$ subject to the relations

$$T_i T_{i+1} T_i = T_{i+1} T_i T_{i+1}, \quad 1 \leq i < n-1$$

$$T_i T_j = T_j T_i, \quad |i - j| \geq 2$$

$$T_0 T_1 T_0 T_1 = T_1 T_0 T_1 T_0$$

$$(T_i - q^{-1})(T_i + q) = 0, \quad 1 \leq i < n-1$$

$$(T_0 - Q^{-1})(T_0 + Q) = 0$$

The relation $T_0 T_1 T_0 T_1 = T_1 T_0 T_1 T_0$ is related to the reflection equation ([Che84, KS92], KSS93 KSS93kss93]). Type B algebras in the context of reflecting boundary conditions have occurred in [BWW18, Kua22]. Alexey Bufetov also has several recent papers on Hecke algebras in probability, such as [Buf20, BC24, BC24b, BB24, BN22]

### 5.2.3 R–matrices and representations

For any vector space $V$, there is a natural representation of the symmetric group $S_n$ on the n–fold tensor product $V^{\otimes n}$ given by

$$\sigma \cdot (v_1 \otimes \cdots \otimes v_n) = v_{\sigma(1)} \otimes \cdots \otimes \cdots v_{\sigma(n)}.$$

In the case of the Hecke algebra $H_q(S_n)$, we need a R–matrix $R: V \otimes V \to V \otimes V$ satisfying

- The braid relations $R_{12} R_{23} R_{12} = R_{23} R_{12} R_{23}$
- The quadratic relation: $(R - q)(R + 1) = 0$ for each i.

In fact, such a $R$ exists (see [Jim86, Jon87, FRT88]). Thus there exists a representation of the Hecke algebra $H_q(S_n)$ on $V^{\otimes n}$ given by having the algebra element $T_i$ act as $R_{i,i+1}$.

Given a solution to a braid relation, generalizing that to a solution to the Yang–Baxter equation is called (Yang)–Baxterization. Generally speaking, this is not an easy task and goes back to e.g. [Jim86, Jon87, ZGB91, Li93].

### 5.2.4 Schur–Weyl Duality

Since the main "body" of this survey is on quantum groups, we should relate Hecke algebras to quantum groups. The classical Schur–Weyl duality is a duality between the action of the symmetric group $S_n$ and the action of the general linear group $GL(V)$ on the $n$-fold tensor product $V^{\otimes n}$. The theorem then states:

#### Theorem

(Schur–Weyl Duality) The actions of $GL(V)$ and $S_n$ on $V^{\otimes n}$ are mutual centralizers, meaning that:

- The centralizer of $S_n$ in $\text{End}(V^{\otimes n})$ is precisely the image of $GL(V)$.
- The centralizer of $GL(V)$ in $End(V^{\otimes n})$ is precisely the image of $S_n$.

The quantum version of Schur–Weyl duality was discovered in [Jim86]:

#### Theorem

(Quantum Schur–Weyl Duality) The actions of $U_q(\mathfrak{gl}(V))$ and $H_q(S_n)$ on $V^{\otimes n}$ are mutual centralizers.

## 5.3 Matrix Product Ansatz

### 5.3.1 Introduction

The Matrix Product Ansatz [DEHP93] is a method for finding stationary measures in an interacting particle system with open boundary conditions. Here, we relate the Matrix Product Ansatz to quantum groups; much of this exposition was motivated by [CRV14].

Let the generator be written in the form

$$B_1 + \sum_{l=1}^{L-1} w_{l,l+1} + \overline{B}_L$$

where we use the subscript notation and $L$ is the number of lattice sites. In other words, $w$ is a local two–site generator, while $B$ acts on the left–most lattice site, and $\overline{B}$ acts on the right–most lattice site.

Denote a particle configuration on $L$ sites by

$$\mathcal{C} = (\tau_1, \ldots, \tau_L),$$

where each $\tau_i$ is a particle variable equal to either 0 or 1, indicating the number of particles at site $i$. The Ansatz takes the form

$$\mathcal{S}(\mathcal{C}) = \frac{1}{Z_L} \langle W | \prod_{i=1}^{L} ((1-\tau_i)E + \tau_i D) | V \rangle$$

where $E$ and $D$ are elements of a non–commutative algebra, and the product is ordered from the index $i = 1$ on the left to $i = L$ on the right.

In [KS97], it was proven that there is a solution of the form

$$w^T\left[\begin{pmatrix}E\\D\end{pmatrix}\otimes\begin{pmatrix}E\\D\end{pmatrix}\right] = \begin{pmatrix}E\\D\end{pmatrix}\otimes\begin{pmatrix}E'\\D'\end{pmatrix} - \begin{pmatrix}E'\\D'\end{pmatrix}\otimes\begin{pmatrix}E\\D\end{pmatrix} \quad (1)$$

and

$$\langle W|B^T\begin{pmatrix}E\\D\end{pmatrix} = \langle W|\left(\begin{pmatrix}E'\\D'\end{pmatrix}\right), \quad \bar{B}^T|\begin{pmatrix}E\\D\end{pmatrix}\rangle = -\begin{pmatrix}E'\\D'\end{pmatrix}|V\rangle \quad (2)$$

where $\overline{E}$ and $D'$ are new non–commuting elements, and the superscript $^T$ denotes transposition. We call these two equations the bulk relations and the boundary relations, respectively. The transposition is there to match terminology between probability papers and mathematical physics papers, since the latter often defines stochastic matrices as having columns summing to 1, rather than rows.

In the case of ASEP, we explicitly take the matrices as

$$B = \begin{pmatrix}-\alpha & \alpha\\ \gamma & -\gamma\end{pmatrix}, \quad w = \begin{pmatrix}0 & 0 & 0 & 0\\ 0 & -q & q & 0\\ 0 & 1 & -1 & 0\\ 0 & 0 & 0 & 0\end{pmatrix}, \quad \bar{B} = \begin{pmatrix}-\delta & \delta\\ \beta & -\beta\end{pmatrix}.$$

The algebraic elements $D$ and $E$ then solve the quadratic relations

$$DE - qED = (1-q)(D+E).$$

This can be expressed in the formulation

$$w^T\left[\begin{pmatrix}E\\D\end{pmatrix}\otimes\begin{pmatrix}E\\D\end{pmatrix}\right] = \begin{pmatrix}E\\D\end{pmatrix}\otimes\begin{pmatrix}1\\-1\end{pmatrix} - \begin{pmatrix}1\\-1\end{pmatrix}\otimes\begin{pmatrix}E\\D\end{pmatrix},$$

so that $\overline{E} = q-1$ and $\overline{D} = 1-q$. Plugging this into the boundary conditions, we obtain

$$\langle W|(\alpha E - \gamma D + q - 1) = 0, \quad (\delta E - \beta D + 1 - q)|V\rangle = 0.$$

In [DEHP93], explicit representations were found, but only for up to three lattice sites.

The extension to arbitrary lattice sizes was done in [San94]. For these purposes, it is useful to define operators $F$ and $F^\dagger$ by (fixing a previous typo in [San94])

$$D = F + 1, \quad E = F^\dagger + 1.$$

Then the quadratic relation between $D$ and $E$ becomes

$$FF^\dagger - qF^\dagger F = 1 - q$$

The operators $F$ and $F^\dagger$ are more commonly known as creation and annihilation operators of the q–deformed harmonic oscillator. There is then an infinite–dimensional representation, with basis $|k\rangle$, going back to [Mac89]. The representation is given by

$$F^\dagger|k\rangle = \{k+1\}_q^{1/2}|k+1\rangle, \quad F|k\rangle = \{k\}_q^{1/2}|k-1\rangle$$

where

$$\{m\}_q = 1 - q^m.$$

The elements $F$ and $F^\dagger$ are adjoint, so

$$\langle k|F = \langle k+1|\{k+1\}_q^{1/2}, \quad \langle k|F^\dagger = \langle k-1|\{k\}_q^{1/2}$$

Next, we need to express $\langle V|$ and $|W\rangle$ in terms of the basis elements. Define the left coefficients $l_k$ and right coefficients $r_k$ by

$$\langle W| = \sum_{k=0}^\infty \langle k|l_k, \quad |V\rangle = \sum_{k=0}^\infty |r_k\rangle.$$

Using that $\langle W|(\alpha E - \gamma D + q - 1) = 0$, we obtain a recurrence relation for the left coefficients

$$l_{-1} = 0, \quad l_0 = 1, \quad \alpha\{k+1\}_q^{1/2}l_{k+1} + (\alpha - \gamma + q - 1)l_k + \gamma\{k\}_q^{1/2}l_{k-1} = 0 \text{ for } k \geq 0$$

Using the adjointness property, then we obtain a similar recurrence relation for the right coefficients

$$r_{-1} = 0, \quad r_0 = 1, \quad \beta\{k+1\}_q^{1/2}r_{k+1} + (\beta - \delta + q - 1)r_k + \delta\{k\}_q^{1/2}r_{k-1} = 0 \text{ for } k \geq 0$$

In principle, solving this recurrence relation gives formulas for the stationary state of ASEP with generic open boundaries on a finite lattice

### 5.3.2 Zamolodchikov Algebra

In this section, we explain how to obtain the bulk relations (1) for the Matrix Product Ansatz from the quantum group. Recall that these are the relations

$$w^T\left[\begin{pmatrix}E\\D\end{pmatrix} \otimes \begin{pmatrix}E\\D\end{pmatrix}\right] = \begin{pmatrix}E\\D\end{pmatrix} \otimes \begin{pmatrix}E'\\D'\end{pmatrix} - \begin{pmatrix}E'\\D'\end{pmatrix} \otimes \begin{pmatrix}E\\D\end{pmatrix},$$

where $w$ is the local generator.

Recall that we assume that the R–matrix $R$ satisfies the usual Yang–Baxter equation

$$R_{1,2}R_{1,3}R_{2,3} = R_{2,3}R_{1,3}R_{1,2}.$$

We will need another property of the R–matrix, called the Markovian property. This states that there exists a vector $v(x)$ such that

$$R^T_{1,2}(zw^{-1})v_1(z)v_2(w) = v_1(z)v_2(w)$$

and the $R$ matrix satisfies the normalization

$$\langle \lambda, \lambda | R^T(x) = \langle \lambda, \lambda |$$

for some highest–weight vector $\langle \lambda, \lambda |$. We explain why this condition is natural.

The universal $R$–matrix of the quantum group has the form

$$R \in U^0 \otimes U^0 + U^{>0} \otimes U^{<0}.$$

The highest weight vector $v$ is defined by $U^{>0}v = 0$, so $v \otimes v$ satisfies the second equality. The first equality essentially amounts to the condition that the $R$–matrix is stochastic, which was discussed in the main set of lecture notes.

The R–matrix needs to be related to the generator $w$. The condition is that if we take the derivative of $R(x)$ with respect to the spectral parameter $x$ and evaluate at 1,

$$PR'(1) = \rho w^T$$

where $P$ is the permutation operator and $\rho$ is some constant. Note that if $R(x)$ is stochastic, then its derivative $R'(x)$ will always have rows summing to 0, from the definition of derivatives.

We will also need the regularity property

$$R(1) = P.$$

As explained in the main body of notes, since the permutation operator solves the braided Yang–Baxter equation, this is a natural condition.

Now, let $T(x)$ be the matrix from the F.R.T construction. Those notes expressed $T$ without a parameter; the parameter comes from the evaluation module in the affine Lie algebra. In this context, we have

$$R_{1,2}(zw^{-1})T_1(z)T_2(w) = T_2(w)T_1(z)R_{1,2}(zw^{-1}).$$

We can now define the Zamolodchikov algebra. Let $v(x) = \begin{pmatrix} v_1(x) \\ v_2(x) \end{pmatrix}$ referencing the vector in the Markovian property. Define

$$A(x) = \begin{pmatrix} \mathcal{X}_1(x) \\ \mathcal{X}_2(x) \end{pmatrix} = T(x)v(x).$$

The Zamolodchikov algebra is the subalgebra of the (dual) quantum group, generated by $\mathcal{X}_1(x), \mathcal{X}_2(x)$.

Applying the RTT relation to $v_1(z)v_2(w)$ we obtain

$$R_{1,2}(zw^{-1})A_1(z)A_2(w) = A_2(w)A_1(z).$$

From here, take the derivative with respect to $z$ and then evaluate at $z = w = 1$. Then we obtain

$$R'_{1,2}(1)A_1(1)A_2(1) + R_{1,2}(1)A'_1(1)A_2(1) = A_2(1)A'_1(1).$$

Plugging in the regularity property $R(1) = P$ and using that $PR'(1) = \rho w$ we then get

$$\rho P w A_1(1)A_2(1) + P A'_1(1)A_2(1) = A_2(1)A'_1(1).$$

Multiplying by $P$ on the left, we obtain

$$w^T A_1(1)A_2(1) = \rho^{-1}(A'_1(1)A_2(1) - A_1(1)A'_2(1))$$

Identifying

$$A_1(1) = \begin{pmatrix} E \\ D \end{pmatrix}, \quad A'_1(1) = \rho \begin{pmatrix} E' \\ D' \end{pmatrix}.$$

we obtain the relation

$$w^T \left[ \begin{pmatrix} E \\ D \end{pmatrix} \otimes \begin{pmatrix} E \\ D \end{pmatrix} \right] = \begin{pmatrix} E \\ D \end{pmatrix} \otimes \begin{pmatrix} E' \\ D' \end{pmatrix} - \begin{pmatrix} E' \\ D' \end{pmatrix} \otimes \begin{pmatrix} E \\ D \end{pmatrix},$$

defining the bulk relations in the Matrix Product Ansatz.

### 5.3.3 Reflection Equation and Ghoshal-Zamolodchikov relations

The reflection equation provides integrable boundary conditions (2) for a spin chain. The reflection equation reads

$$R_{1,2}(zw^{-1})K_1(z)R_{2,1}(zw)K_2(w) = K_2(w)R_{1,2}(zw)K_1(z)R_{2,1}(zw^{-1})$$

and the boundary conditions are

$$\langle W|(BA(1) - \rho^{-1}A'(1)) = 0, \quad (\overline{B}A(1) - \rho^{-1}A'(1))|V\rangle = 0.$$

A visual representation is given in [CRV14]:

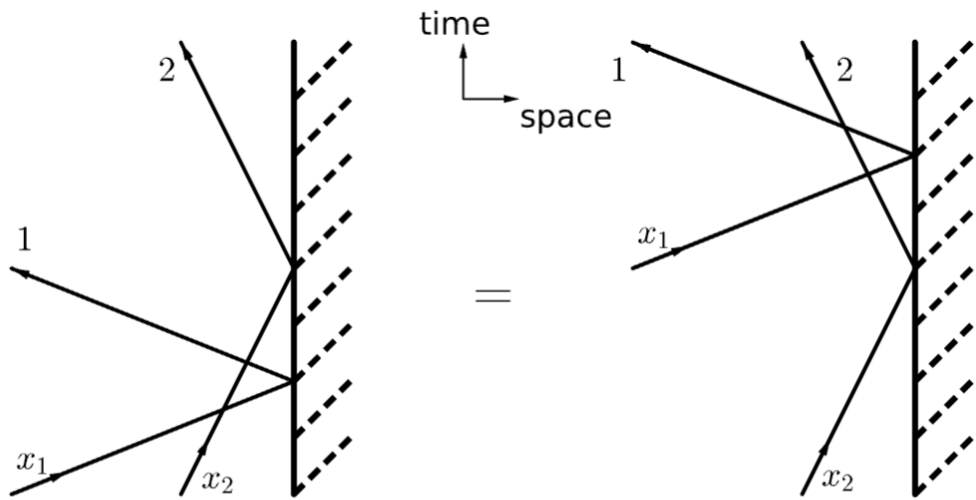

Two rays reflecting off a barrier, with one-dimensional space displayed in the horizontal direction and time in the vertical direction. The rays are labeled 1 and 2. Ray 1 moves faster than ray 2. On the left side of the image, ray 1 crosses ray 2, then bounces off the wall, then crosses ray 2 again, and then ray 2 reflects off the wall. On the right side of the image, ray 2 bounces off the wall, then crosses ray 1, then ray 1 reflects off the wall, then ray 1 crosses ray 2 again. The two sides are equal.

*Figure 7: The reflection equation. If the crossings are read from bottom to top, one obtains the reflection equation.*

Given an R-matrix, solutions of these equations are called reflection matrices and denoted $K(z)$. In general, these solutions are not unique. In our context, we will consider two solutions denoted $K(z)$ and $\overline{K}(z)$, for the left and right boundaries.

We make a few assumptions about the reflection matrix $K(z)$. First, we assume the regularity condition

$$K(1) = Id.$$

This is a natural condition, as it corresponds to closed boundary conditions. Additionally, K should be related to the local generators at the boundaries:

$$K'(1) = 2\rho B, \quad \overline{K}'(1) = -2\rho \overline{B}$$

The Ghoshal-Zamolodchikov relations, which we will abbreviate as the GZ relations, are given by

$$\langle \Omega | (K(z)A(z^{-1}) - A(z)) = 0, \quad (\overline{K}(z)A(z^{-1}) - A(z))|\Omega\rangle = 0.$$

To explain the connection to the reflection equation, we calculate $\langle \Omega | A_2(z_2) A_1(z_1)$ in two different ways. Assuming the GZ–relations, we first have

$$\langle \Omega | A_2(z_2) A_1(z_1) = R_{1,2}(z_1 z_2^{-1}) \langle \Omega | A_1(z_1) A_2(z_2)$$

$$= R_{1,2}(z_1 z_2^{-1}) K_1(z_1) \langle \Omega | A_1(z_1^{-1}) A_2(z_2)$$

$$= R_{1,2}(z_1 z_2^{-1}) K_1(z_1) R_{2,1}(z_1 z_2) \langle \Omega | A_2(z_2) A_1(z_1^{-1})$$

$$= R_{1,2}(z_1 z_2^{-1}) K_1(z_1) R_{2,1}(z_1 z_2) K_2(z_2) \langle \Omega | A_2(z_2^{-1}) A_1(z_1^{-1})$$

By a similar sequence of calculations,

$$\langle \Omega | A_2(z_2) A_1(z_1) = K_2(z_2) R_{1,2}(z_1 z_2) K_1(z_1) R_{2,1}(z_1 z_2^{-1}) \langle \Omega | A_2(z_2^{-1}) A_1(z_1^{-1}).$$

Thus, the reflection equation is a sufficient condition to guarantee the consistency of the GZ relations.

If we take the derivative of the GZ–relations at $z = 1$ we obtain

$$\langle \Omega | (K'(1)A(1) - K(1)A'(1) - A'(1)) = 0 \quad ((\overline{K}'(1)A(1) - \overline{K}(1)A'(1) - A'(1))|\Omega\rangle = 0.$$

Defining $\langle W | = \langle \Omega |$ and $|V\rangle = |\Omega\rangle$, we obtain

$$\langle W | (BA(1) - \rho^{-1} A'(1)) = 0, \quad (\overline{B}A(1) - \rho^{-1} A'(1))|V\rangle = 0.$$

Thus, we obtain the boundary conditions in the Matrix Product Ansatz.

### 5.4.4 Explicit Formulas

In this section, we elaborate on explicit formulas for the spectral–dependent R–matrices and their relation to ASEP. We define the matrix entries

$$T(z) = \begin{pmatrix} T_{11}(z) & T_{12}(z) \\ T_{21}(z) & T_{22}(z) \end{pmatrix}.$$

We write these entries as a power series

$$T_{ij}(z) = \sum_{n=0}^{\infty} T_{ij}^{(n)} z^n.$$

We will consider certain evaluation maps, which allow for the infinite–dimensional affine Lie algebra to be evaluated on finite–dimensional modules. We can truncate the infinite series:

$$\begin{pmatrix} T_{11}^{(0)} & 0 \\ T_{21}^{(0)} & 1 \end{pmatrix} + z \begin{pmatrix} 1 & T_{12}^{(0)} \\ 0 & T_{22}^{(0)} \end{pmatrix}$$

For context on this evaluation map, one can refer to e.g. [ZG93], in particular the lemma on the evaluation module.

For ASEP, the R–matrix can be expressed as

$$\begin{pmatrix} 1 & 0 & 0 & 0 \\ 0 & \dfrac{(z-1)q}{qz-1} & \dfrac{q-1}{qz-1} & 0 \\ 0 & \dfrac{(q-1)z}{qz-1} & \dfrac{z-1}{qz-1} & 0 \\ 0 & 0 & 0 & 1 \end{pmatrix}$$

Two solutions of the reflection equation are:

$$K(x) = \begin{pmatrix} \dfrac{(-x\alpha + x\gamma + q + \alpha - \gamma - 1)x}{x^2\gamma + qx + x\alpha - x\gamma - x - \alpha} & \dfrac{(x^2-1)\gamma}{x^2\gamma + qx + x\alpha - x\gamma - x - \alpha} \\ \dfrac{\alpha(x^2-1)}{x^2\gamma + qx + x\alpha - x\gamma - x - \alpha} & -\dfrac{(-qx - x\alpha + x\gamma + x + \alpha - \gamma)}{x^2\gamma + qx + x\alpha - x\gamma - x - \alpha} \end{pmatrix},$$

and

$$\overline{K}(x) = \begin{pmatrix} \dfrac{(x\delta - x\beta + q - \delta + \beta - 1)x}{-x^2\beta + qx - x\delta + x\beta - x + \delta} & -\dfrac{(x^2-1)\beta}{-x^2\beta + qx - x\delta + x\beta - x + \delta} \\ -\dfrac{(x^2-1)\delta}{-x^2\beta + qx - x\delta + x\beta - x + \delta} & \dfrac{qx - x\delta + x\beta - x + \delta - \beta}{-x^2\beta + qx - x\delta + x\beta - x + \delta} \end{pmatrix}.$$

We also have $\rho = (q-1)^{-1}$.

The Markovian property is satisfied by

$$v(z) = f(z)\left(\dfrac{ax}{b}\right)$$

where $f(z)$ is an arbitrary function, and $a, b$ are arbitrary constants. Using the evalution map, the generators of the Zamolodchikov algebra are

$$A(x) = f(x) \begin{pmatrix} ax^2 + x(aT_{11}^{(0)} + bT_{12}^{(1)}) \\ x(aT_{21}^{(0)} + bT_{22}^{(1)}) + b \end{pmatrix}.$$

Choosing $a = b = 1$ and $f(x) = \dfrac{1}{x}$, one gets

$$E = 1 + T_{11}^{(0)} + T_{12}^{(1)}, \quad \overline{E} = q - 1, \quad D = T_{21}^{(0)} + T_{22}^{(1)} + 1 \quad \overline{D} = 1 - q$$

## 5.4 Orthogonal Polynomials

### 5.4.1 Introductory example

While this appendix does depend on the main document, I will also try to make this as self–contained as possible. To that end, I will provide a well-known motivation coming from quantum harmonic oscillators. Even from this most basic example, we will be able to conclude that **orthogonal polynomials occur in the matrix entries of representations.** Since vertex weights are entries of the $R$–matrix coming from representation theory, we expect that vertex weights are orthogonal polynomials.

In contrast to classical mechanics where position and momentum are either vectors or scalars, in quantum mechanics these are operators. In general, these operators do not commute, and the non–commutativity of the position and momentum operators is related to Heisenberg's uncertainty principle. If the two operators commuted, then both position and momentum could be measured simultaneously. Furthermore, states in a quantum mechanical system are actually "wavefunctions," which are modeled as an element of some Hilbert space (on which the operators act).

In one dimension, for example, the position operator $\hat{x}$ is multiplication by $x$ while the momentum operator $\hat{p}$ is $-i\hbar \frac{\partial}{\partial x}$. One can verify that if $f_n$ denotes the function $x^n$ then

$$\hat{x}\hat{p} f_n = -i\hbar n f_n, \quad \hat{p}\hat{x} = -i\hbar(n+1) f_n,$$

and so (using Lie bracket notation)

$$[\hat{p}, \hat{x}] f_n = -i\hbar f_n \neq 0.$$

More generally, the time–dependent Schrödinger's equation states that for the Hamiltonian operator $\hat{H}$ and a time-dependent state $\Psi(t)$ we have

$$i\hbar \frac{d}{dt} |\Psi(t)\rangle = \hat{H} |\Psi(t)\rangle$$

where $i$ is the imaginary unit and $\hbar$ is Planck's constant. In the time–independent version, we consider the case when $\Psi(t)$ is a stationary state, so $\Psi(t) = \Psi$ for all times $t$. Then Schrödinger's equation tells us that

$$\hat{H} |\Psi\rangle = E |\Psi\rangle,$$

where is the energy of the state. Note that this means that the energy can only take discrete, "quantized" values. Generally speaking, the Hamiltonian can be expressed as

$$\hat{H} = \hat{T} + \hat{V}$$

where $\hat{T}$ and $\hat{V}$ are the kinetic and potential energy.

### 5.4.2 Quantum Harmonic Oscillator and Hermite Polynomials

For the harmonic oscillator, we get (in analogy with the classical harmonic oscillator)

$$\frac{\hat{p}^2}{2m} + \frac{1}{2}k\hat{x}^2.$$

Plugging in $\hat{p}, \hat{x}$ we get a differential equation which can be explicitly solved. More specifically, the steady states are

$$\psi_n(x) = \frac{1}{\sqrt{2^n n!}} \left(\frac{m\omega}{\pi\hbar}\right)^{1/4} e^{-\frac{m\omega x^2}{2\hbar}} H_n\left(\sqrt{\frac{m\omega}{\hbar}} x\right), \qquad n = 0,1,2,\dots$$

where $H_n$ are the physicists' Hermite polynomials

$$H_n(z) = (-1)^n e^{z^2} \frac{d^n}{dz^n}\left(e^{-z^2}\right).$$

The corresponding energy levels are:

$$E_n = \hbar\omega\left(n + \frac{1}{2}\right).$$

The Hermite polynomials are orthogonal polynomials in the sense that

$$\int_{-\infty}^{\infty} H_m(x)H_n(x)e^{-x^2}dx = \sqrt{\pi}2^n n!\, \delta_{nm}.$$

They satisfy the generating function

$$e^{2xt-t^2} = \sum_{n=0}^{\infty} H_n(x)\frac{t^n}{n!}$$

and the recurrence relation

$$H_{n+1}(x) = 2xH_n(x) - 2nH_{n-1}(x).$$

### 5.4.3 CCR relations and Lie algebras

Now, what does this have to do with Lie algebras? To answer this, we approach the quantum harmonic oscillator using "ladder operators." Define two abstract algebra elements $\hat{a}$ and $\hat{a}^\dagger$ with canonical commutation relation

$$[\hat{a}, \hat{a}^\dagger] = 1.$$

There is a representation onto an infinite–dimensional vector space with bases $|n\rangle, n = 0,1,2$ given by

$$\hat{a}^\dagger|n\rangle = \sqrt{n+1}|n+1\rangle$$
$$\hat{a}|n\rangle = \sqrt{n}|n-1\rangle.$$

For this reason, $\hat{a}$ and $\hat{a}^\dagger$ are called lowering and raising operators, respectively. We also introduce a "number operator"

$$N = \hat{a}^\dagger \hat{a}$$
$$N|n\rangle = n|n\rangle.$$

In the context of the quantum harmonic oscillator, these algebra elements can be represented as

$$\hat{a} = \sqrt{\frac{m\omega}{2\hbar}} \left( \hat{x} + \frac{i}{m\omega} \hat{p} \right)$$

$$\hat{a}^\dagger = \sqrt{\frac{m\omega}{2\hbar}} \left( \hat{x} - \frac{i}{m\omega} \hat{p} \right)$$

Then

$$\hat{H} = \hbar\omega \left( N + \frac{1}{2} \right),$$

so the energy levels are $E_n = n + \frac{1}{2}$, as we found before.

To relate the two methods in more mathematical formalism, we define the Hilbert space $L^2(\mathbb{R}, w)$ of square–integrable functions on $\mathbb{R}$ with respect to the weight $w(x) = e^{-x^2/2}$. Then this Hilbert space has an orthogonal basis $\{H_n : n \geq 0\}$ consisting of the Hermite polynomials. If we let $\mathcal{A}$ denote the algebra generated by $\hat{a}, \hat{a}^\dagger$ then we just defined a representation of $\mathcal{A}$ on $L^2(\mathbb{R}, w)$ where $|n\rangle = H_n$. Note that we can define a map from $\mathcal{A}$ to $B(L^2(\mathbb{R}, w))$. Using the generating function and recurrence relation, we have

$$xH_n(x) = \frac{1}{2} H_{n+1}(x) + nH_{n-1}(x)$$

and

$$\frac{d}{dx} H_n(x) = 2nH_{n-1}(x).$$

Then the ladder operators act as

$$\hat{a} = \frac{1}{\sqrt{2}} \left( x + \frac{d}{dx} \right), \quad \hat{a}^\dagger = \frac{1}{\sqrt{2}} \left( x - \frac{d}{dx} \right).$$

Then using the formulas above, these act on the wavefunctions as

$$\hat{a}^\dagger \psi_n = \sqrt{n+1}\, \psi_{n+1}$$
$$\hat{a} \psi_n = \sqrt{n}\, \psi_{n-1}$$

Next, we relate Lie algebras to ladder operators using the Jordan–Schwinger realization. It may be more illustrative to demonstrate this first in the case of $\mathfrak{sl}_2$. Let $\mathcal{A}_n$ be generated by $\hat{a}_i, \hat{a}_i^\dagger, i = 1,2,\ldots,n$ with the commutation relations

$$[\hat{a}_i, \hat{a}_j^\dagger] = 1\delta_{ij}, \quad [\hat{a}_i, \hat{a}_j] = [\hat{a}_i^\dagger, \hat{a}_j^\dagger] = 0 \text{ for } i \neq j.$$

The Jordan–Schwinger realization of $\mathfrak{sl}_2$ is then an injection from $\mathfrak{sl}_2$ to $\mathcal{A}_2$. In a typical physics introduction, this example is discussed in the context of angular momentum. In this context, we recall the Pauli spin matrices

$$\sigma^1 = \begin{pmatrix} 0 & 1 \\ 1 & 0 \end{pmatrix}, \qquad \sigma^2 = \begin{pmatrix} 0 & -i \\ i & 0 \end{pmatrix}, \qquad \sigma^3 = \begin{pmatrix} 1 & 0 \\ 0 & -1 \end{pmatrix}.$$

which span $\mathfrak{sl}_2$. The change of basis to $e, f, h$ is given by

$$e = \frac{\sigma_1 + i\sigma_2}{2}$$

$$\sigma_1 = e + f$$

$$f = \frac{\sigma_1 - i\sigma_2}{2}$$

$$\sigma_2 = -i(e - f)$$

$$h = \sigma_3$$

$$\sigma_3 = h$$

The Pauli matrices are then written in terms of the ladder operators as

$$\sigma_1 = \hat{a}_1 \hat{a}_2^\dagger + \hat{a}_1^\dagger \hat{a}_2, \quad \sigma_2 = -i(\hat{a}_1^\dagger \hat{a}_2 - \hat{a}_1 \hat{a}_2^\dagger), \quad \sigma_3 = \hat{a}_1^\dagger \hat{a}_1 - \hat{a}_2^\dagger \hat{a}_2$$

By comparing the last formula to the change of basis, it is reasonable to conjecture that the elements $\hat{a}_1^\dagger \hat{a}_2, \hat{a}_1 \hat{a}_2^\dagger, \hat{a}_1^\dagger \hat{a}_1 - \hat{a}_2^\dagger \hat{a}_2$ then satisfy the same relations as $e, f, h \in \mathfrak{sl}_2$, respectively. Indeed this is the case; for example, the relation $[h, e] = 2e$ can be seen in

$$(\hat{a}_1^\dagger \hat{a}_1 - \hat{a}_2^\dagger \hat{a}_2) \cdot (\hat{a}_1^\dagger \hat{a}_2) = \hat{a}_1^\dagger \hat{a}_1 \hat{a}_1^\dagger \hat{a}_2 - \hat{a}_2^\dagger \hat{a}_2 \hat{a}_1^\dagger \hat{a}_2$$

$$= \hat{a}_1^\dagger (\hat{a}_1^\dagger \hat{a}_1 + 1) \hat{a}_2 - \hat{a}_1^\dagger (\hat{a}_2 \hat{a}_2^\dagger - 1) \hat{a}_2$$

$$= \hat{a}_1^\dagger \hat{a}_2 \hat{a}_1^\dagger \hat{a}_1 + \hat{a}_1^\dagger \hat{a}_2 - \hat{a}_1^\dagger \hat{a}_2 \hat{a}_2^\dagger \hat{a}_2 + \hat{a}_1^\dagger \hat{a}_2$$

$$= (\hat{a}_1^\dagger \hat{a}_2) \cdot (\hat{a}_1^\dagger \hat{a}_1 - \hat{a}_2^\dagger \hat{a}_2 + 2).$$

More generally, there is an injection from any finite–dimensional Lie algbra $\mathfrak{g}$ to some $\mathcal{A}_n$. Since there is an injection $ad: \mathfrak{g} \to \mathfrak{gl}(\mathfrak{g})$, without loss of generality we can take $\mathfrak{g} = \mathfrak{gl}(n)$. The Jordan–Schwinger realization is then an injection of Lie algebras defined by

$$M \to \sum_{i,j=1}^{n} \hat{a}_i^\dagger M_{ij} \hat{a}_j.$$

One can check that this generalizes the case when $\mathfrak{g} = \mathfrak{sl}_2$.

As a result, if there is a Jordan–Schwinger realization $\mathfrak{g} \hookrightarrow \mathcal{A}_n$ then any representation of $\mathcal{A}_n$ pulls back to a representation of $\mathfrak{g}$. Since we have already seen the Hermite polynomials appear as an orthonormal basis of a representation of $\mathcal{A}_1$, it is not unreasonable to expect orthogonal polynomials to appear in the representations of $\mathfrak{g}$. Note that if we had instead picked a basis of $1, x, x^2, \ldots$, then after a change of basis, we would see the Hermite **polynomials in the matrix elements of representations** of $\mathcal{A}_1$.

Of course, the Jordan–Schwinger realization does not produce interesting examples of Lie algebras inside $\mathcal{A}_1$, so we will have to go to higher–dimensional examples. This will result in generalizations of the Hermite polynomials, which leads into the next section.

### 5.4.4 Hierarchy of Orthogonal Polynomials

The Askey–Wilson hierarchy is a way of organizing orthogonal polynomials of hypergeometric or basic hypergeometric type into a hierarchy. Below is the hierarchy for basic hypergeometric type (courtesy of Koornwinder):

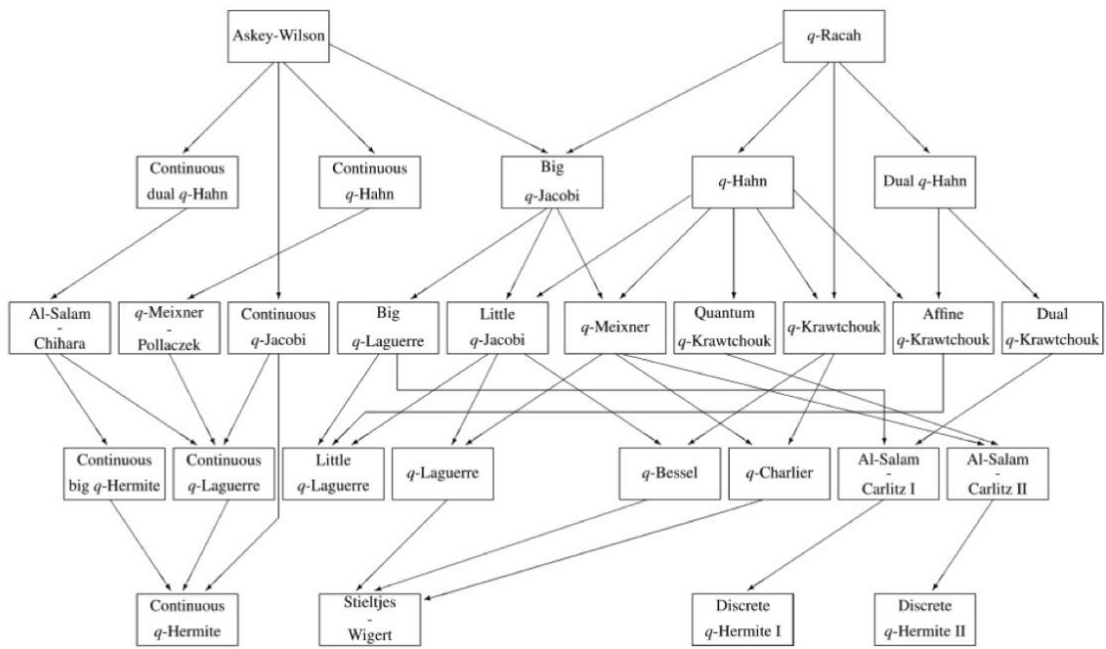

Show Alt Text

There is a more accessible version of this flowchart of the Askey-Wilson hierarchy , written as a webpage using ARIA flow-to.

*Figure 8: The Askey-Wilson hierarchy.*

Note that Hermite polynomials are in this flowchart, as a $q \to 1$ limit of the various $q$–Hermite polynomials. We note that the $q$–deformation here is a different quantization than the quantization in the quantum harmonic oscillator; roughly, this corresponds to quantizing the product in one case, and quantizing the co–product in the other case.

### 5.4.5 Definitions

Define the $q$–Pochhammer symbol by

$$(a;q)_n = (1-a)(1-aq)\cdots(1-aq^{n-1})$$

and

$$(a_1,\ldots,a_r)_n = (a_1;q)_n \cdots (a_r;q)_n$$

The basic hypergeometric series, or $q$–hypergeometric series is defined by ([KLSBook])

$$_r\phi_s\begin{bmatrix} a_1 & ;a_2 & ;\ldots & ;a_j \\ b_1 & ;b_2 & ;\ldots & ;b_k \end{bmatrix};q,z\Big] := \sum_{k=0}^{\infty} \frac{(a_1,\ldots,a_r;q)_k}{(b_1,\ldots,b_s;q)_k}(-1)^{(1+s-r)k}q^{(1+s-r)\binom{k}{2}}\frac{z^k}{(q;q)_k}$$

The various orthogonal polynomials are defined in terms of the $q$–hypergeometric series. Note that if $a_1 = q^{-n}$ then the series terminates, justifiably calling it a polynomial. One example are $q$–Racah polynomials, which are defined by

$$R_n(\mu(x);\alpha,\beta,\gamma,\delta|q) = {}_4\varphi_3\left(\begin{matrix} q^{-n},\alpha \\ \beta q^{n+1},q^{-x},\gamma\delta q^{x+1}\alpha q, \beta\delta q, \gamma q \end{matrix};q,q\right), n = 0,1,2,\ldots,N.$$

and

$$\alpha q = q^{-N} \quad \text{or} \quad \beta\delta q = q^{-N} \quad \text{or} \quad \gamma q = q^{-N}.$$

Note that $_4\varphi_3$ depends on seven variables whereas the $q$–Racah depends on six variables; if the relation $a_1 a_2 a_3 a_4 = q^{-1} b_1 b_2 b_3$ holds, and one of $\{b_1, b_2, b_3\}$ is of the form $q^{-N}$, then the $_4\varphi_3$ function can be written as a $q$–Racah.

It turns out that the weights of the $R$–matrix of $\mathcal{U}_q(\widehat{sl_2})$ are given by the $q$–Racah polynomials. Before stating that result, we define some notation. We will let $m$ denote the maximum number of particles that can enter from the vertical direction, and let $l$ denote the maximum number of particles that can enter from the horizontal direction. So the stochastic matrix is an operator on $V_l \otimes V_m$ where $V_l$ is $l+1$–dimensional and $V_m$ is $m+1$–dimensional. We denote a matrix entry by

$$[S(z)]_{j_1,k_1}^{j_2,k_2}, \quad 1 \leq j_1, j_2 \leq l, \quad 1 \leq k_1, k_2 \leq m.$$

This corresponds to this weight:

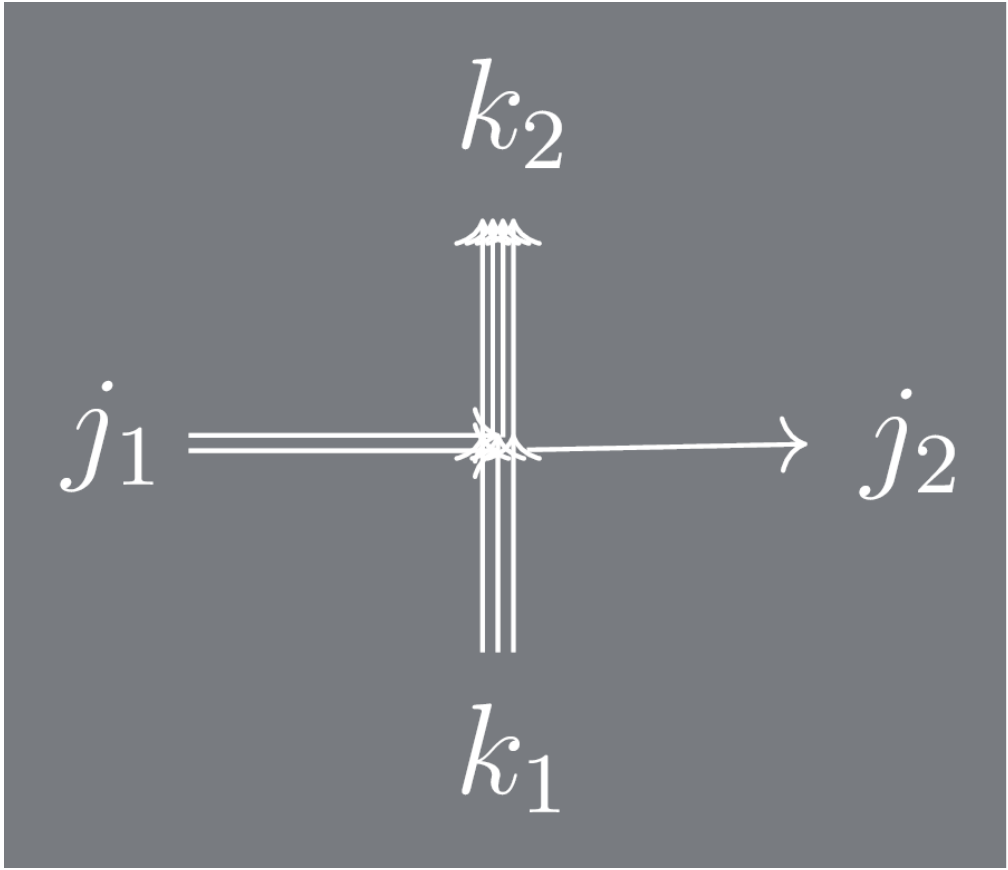

A vertex weight labeled by $j_1$ at the left and $j_2$ at the right and $k_1$ at the bottom and $k_2$ at the top. The number of arrows coming in from the left and bottom equals the number of arrows going out the top and right.

*Figure 9: The vertex weight corresponding to $[S(z)]_{j_1,k_1}^{j_2,k_2}$*

Of course, the weights of the $R$–matrix are generally not stochastic. There is a type of transformation that will make it stochastic.

### Definition

Let $G_{lm}: V_l \otimes V_m \to V_l \otimes V_m$ be a diagonal map (with respect to some choice of basis). Let $P: V_l \otimes V_m \to V_m \otimes V_l$ be the permutation map $P(v \otimes w) = w \otimes v$. Then a gauge transformation of a matrix $R: V_l \otimes V_m \to V_l \otimes V_m$ is another matrix $S: V_l \otimes V_m \to V_l \otimes V_m$ such that

$$S = PG_{ml}^{-1}PRG_{lm}$$

We will state a theorem from [CP16] and [Kua18]; see also [KMMO16] and [BM16] for higher rank generalizations; and the previous papers of [Pov13] and [Man14].

### Theorem

There exists a gauge transformation of the $R$–matrix of $U_q(\widehat{sl}_2)$ such that the stochastic vertex weights are given by (where $v$ depends on $m$ as $v = q^{-m}$):

$$\left[S^{l,m}(z)\right]_{j_1,k_1}^{j_2,k_2} = \mathbf{1}_{j_1+k_1=j_2+k_2} q^{\frac{2k_1-k_1^2}{4} - \frac{2k_2-k_2^2}{4} + \frac{j_2^2+j_1^2}{4} + \frac{j_2(k_2-1)+j_1k_1}{2}}$$

$$\times \frac{v^{k_1-j_2}\alpha^{k_2-k_1+j_2}(-\alpha v^{-1};q)_{k_2-j_1}}{(q;q)_{j_2}(-\alpha;q)_{j_2+k_2}(\beta\alpha^{-1}q^{1-k_1};q)_{k_1-k_2}} R_{j_1}\left(\mu(j_2); \middle| q^{k_2-j_1}, -vq^{k_2}z, zq^lq^{-1-j_2-k_2}, -q^{k_2-1}z\right)$$

$$\times (v;q)_{j_1}\left(q^{1+k_2-j_1};q\right)_{j_1}\left(q^{l+1-j_2-k_2};q\right)_{j_1}$$

Note that the reason I used $v = q^{-m}$ is that the value $v$ can be extended to arbitrary complex numbers through analytic continuation. I believe this can be done algebraically through Verma modules, but I'm not sure.

The proof of this theorem is essentially in two parts. In Theorem 3.15 of [CP16], the authors prove that applying Rogers–Pitman intertwining results in stochastic vertex weights expressed in terms of the $q$–Racah polynomials. In [Kua18], I prove that applying Rogers–Pitman intertwining "commutes" with the gauge transformation, in the sense that first applying the gauge transformation and then applying fusion results in the same matrix as first applying fusion and then applying the gauge transformation. Thus, the stochastic matrix $S$ has the necessary algebraic properties to prove Markov duality.

I'll only sketch the argument. The proof that Rogers–Pitman intertwining "commutes" with the gauge transformation is in Theorem 3.4 of [Kua18], and is a long calculation that is not particularly enlightening. The proof that the $q$–Racah polynomials occur is due to a recurrence relation that comes from fusion. Before that, we consider the $l = 1$ case, so that the weights appear in Figure 10.

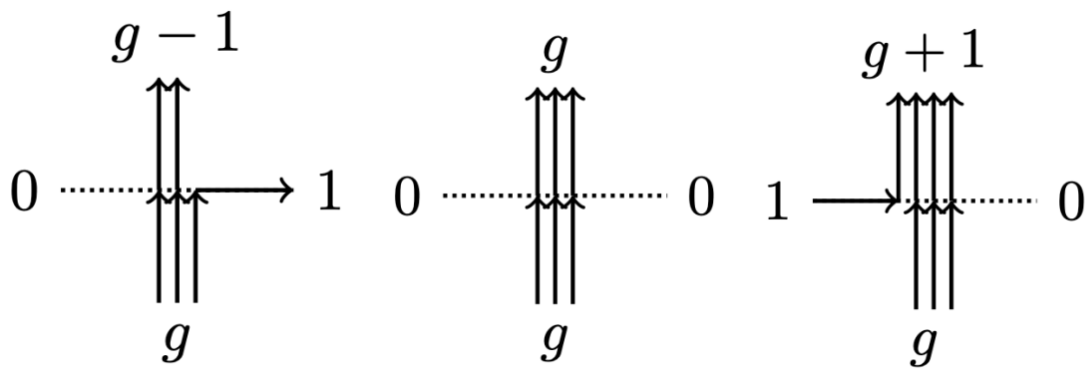

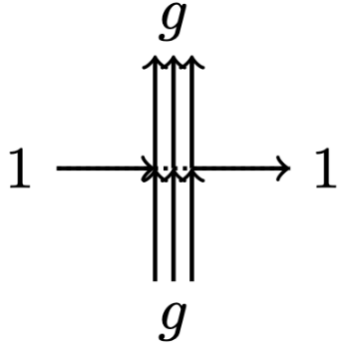

Four vertices; the top right vertex has g arrows coming from the bottom and g arrows going out the top. The top left vertex has g arrows coming from the bottom, and one arrow going out the right and $g - 1$ arrows coming out the top. The bottom right vertex has g arrows coming from the bottom and 1 arrow coming from the left, and g arrows going out the right and 1 arrow doing out the top. The bottom left vertex has g arrows coming from the bottom and 1 arrow coming from the left, and $g + 1$ arrows coming out the top.

Figure 10: The top right vertex has weight $\frac{(q^{m+1}-q^{2g}z)}{q^{m+1}-z}$ and the top left vertex has weight $\frac{z(q^{2g}-1)}{q^{m+1}-z}$. The bottom right vertex has weight $\frac{q^{2g-m+1}-z}{q^{m+1}-z}$ and the bottom left vertex has weight $\frac{q^{m+1}-q^{2g-m+1}}{q^{m+1}-z}$

Those weights are simple enough to calculate. Next, we fuse in the vertical direction. Figure 11 below shows the intuition for the argument. If $j_1$ arrows come in from the left, where at most $l$ arrows are allowed, we choose random places for the arrows across two vertices where at most $l - 1$ and 1 arrows are allowed. The random placement has a nice combinatorial description. Define the $q$–deformed integer, factorial and binomial by

$$[n]_q = \frac{1-q^n}{1-q}, \quad [n]_q^! = [1]_q \cdots [n]_q, \quad \binom{l}{j}_q = \frac{[l]_q^!}{[j]_q^! [l-j]_q^!}.$$

Then we have the relation

$$\binom{l}{j}_q = q^j \binom{l-1}{j}_q + \binom{l-1}{j-1}_q.$$

So then

$$\mathbf{P}(0) = \frac{q^j \binom{l-1}{j}_q}{\binom{l}{j}_q}, \quad \mathbf{P}(1) = \frac{\binom{l-1}{j-1}_q}{\binom{l}{j}_q},$$

which simplifies to the same probability in Lemma 3.13 of [CP16].

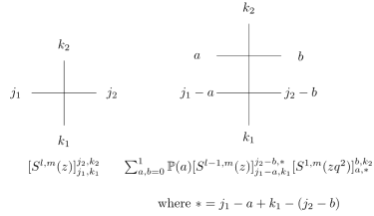

A vertex labeled with $k_1$ at the bottom, $k_2$ at the top, $j_1$ at the left and $j_2$ at the right. Underneath that vertex is $[S^{l,m}(z)]_{j_1,k_1}^{j_2,k_2}$. Next to that vertex is two connected vertices shaped like ⊢, again with $k_1$ at the bottom. Above that is a row with $j_1 - a$ on the left and $j_2 - b$ on the right. Above that is another row with a on the left and b on the right. Again there is $k_2$ at the top. Below that vertex is $\sum_{a,b=0}^{1} \mathbb{P}(a)[S^{l-1,m}(z)]_{j_1-a,k_1}^{j_2-b,*}[S^{1,m}(zq^2)]_{a,*}^{b,k_2}$ where $* = j_1 - a + k_1 - (j_2 - b)$.

Figure 11: *A visual heuristic of fusion of vertex weights.*

In this way, you obtain a recurrence relation in the parameter $l$. Lemma 3.13 of [CP16] shows that the $q$–Racah formula satisfies the same recurrence relation.

# References

**[ADHR94]**
F. C. Alcaraz, M. Droz, M. Henkel, and V. Rittenberg, Ann. Phys. (N.Y.)230:250 (1994).

**[BDJ99]**
Baik, J.; Deift, P.; Johansson, K. (1999), "On the distribution of the length of the longest increasing subsequence of random permutations", Journal of the American Mathematical Society, 12 (4): 1119–1178.

**[Bax72]**
R. J. Baxter. One-dimensional anisotropic Heisenberg chain. Ann. Physics, 70:323–337, 1972.

**[BWW18]**
Huanchen Bao, Weiqiang Wang, and Hideya Watanabe, "Multiparameter quantum Schur duality of type B," Proceedings of the American Mathematical Society, Vol. 146 (2018), no. 9, pp. 3203-3216

**[BM16]**
Gary Bosnjak, Vladimir V. Mangazeev, Construction of R-matrices for symmetric tensor representations related to $U_q(\widehat{sl_n})$. J. Phys. A: Math. Theor. 49 (2016) 495204

**[Buf20]**
A. Bufetov, *Interacting particle systems and random walks on Hecke algebras*, arXiv:2003.02730 (2020).

**[BC24]**

A. Bufetov and K. Chen, *Mallows product measure*, Electron. J. Probab. **29** (2024), 1–33, arXiv:2402.09892.

**[BC24b]**

A. Bufetov and K. Chen, *Local central limit theorem for Mallows measure*, arXiv:2409.10415 (2024).

**[BB24]**

A. Borodin and A. Bufetov, *ASEP via Mallows coloring*, arXiv:2408.16585 (2024).

**[BN22]**

A. Bufetov and P. Nejjar, *Cutoff profile for ASEP on a segment*, Probab. Theory Related Fields **183** (2022), 229–253, arXiv:2012.14924.

**[CGRS16]**

Carinci, G., Giardiná, C., Redig, F. Redig, T. Sasamoto. A generalized asymmetric exclusion process with stochastic duality. Probab. Theory Relat. Fields 166, 887–933 (2016).

**[CGX91]**

Y. Cheng, M.L. Ge and K. Xue, Yang–Baxterization of Braid Group Representations, Commun. Math. Phys. 136 (1991) 195.

**[Che84]**

I. Cherednik, *Factorizing particles on a half-line and root systems*, Theoret. and Math. Phys. **61** (1984), 977–983.

**[CP16]**

Corwin, I., Petrov, L. Stochastic Higher Spin Vertex Models on the Line. Commun. Math. Phys. 343, 651–700 (2016).

**[Cox35]**

H. S. M. Coxeter, The complete enumeration of finite groups of the form $r_i^2 = (r_i r_j)^{k_{ij}} = 1$, J. London Math. Soc. **10** (1935), 21–25.

**[CRV14]**

N Crampe, E Ragoucy and M Vanicat, "Integrable approach to simple exclusion processes with boundaries. Review and progress." Journal of Statistical Mechanics: Theory and Experiment, Volume 2014, November 2014

**[DEHP93]**

B Derrida, M R Evans, V Hakim and V Pasquier Exact solution of a 1D asymmetric exclusion model using a matrix formulation Journal of Physics A: Mathematical and General, Volume 26, Number 7, 1993.

**[Dri85]**
Drinfeld, V. G. (1985), "Hopf algebras and the quantum Yang–Baxter equation", Doklady Akademii Nauk SSSR, 283 (5): 1060–1064,

**[Dri87]**
Congress of Mathematicians, Vol. 1, 2 (Berkeley, Calif., 1986), 798–820, Amer. Math. Soc., Providence, RI, 1987.

**[FRT88]**
L. D. Faddeev, N. Yu. Reshetikhin, L. A. Takhtajan, Quantization of Lie groups and Lie algebras, Algebraic analysis, Vol. I, Academic Press, Boston, MA, (1988), 129-139.

**[GS92]**
L.-H. Gwa and H. Spohn, Phys. Rev. A 46:844 (1992).

**[Hei28]**
Heisenberg, W. (1 September 1928). "Zur Theorie des Ferromagnetismus" [On the theory of ferromagnetism]. Zeitschrift für Physik (in German). 49 (9): 619–636.

**[Hop41]**
Hopf, Heinz (1941). Über die Topologie der Gruppen–Mannigfaltigkeiten und ihre Verallgemeinerungen. Ann. of Math. 2 (in German). 42 (1): 22–52.

**[Jim85]**
Jimbo, Michio (1985), "A q-difference analogue of $U(g)$ and the Yang-Baxter equation", Letters in Mathematical Physics, 10 (1): 63–69,

**[Jim86]**
M. Jimbo, *A q-analogue of $U(\mathfrak{gl}(N+1))$, Hecke algebra, and the Yang-Baxter equation*, Lett. Math. Phys. **11** (1986), 247–252.

**[Joh00]**
Johansson, K., "Shape fluctuations and random matrices", Communications in Mathematical Physics, 209 (2): 437–476, 2000.

**[Jon87]**
V. F. R. Jones, Hecke Algebra Representations of Braid Groups and Link Polynomials, Annals of Mathematics, Vol. 126, No. 2 (Sep., 1987), pp. 335-388 (54 pages)

**[Jon90]**
V.F.R. Jones, "Baxterisation", **Int. J. Mod. Phys**. B 4 (1990) 701, proceedings of "Yang–Baxter equations, conformal invariance and integrability in statistical mechanics and field theory", Canberra, 1989.

**[KS60]**
J.G. Kemeny and J.L. Snell, Finite Markov Chains, D. Van Nostrand Company, Princeton, New Jersey, 1960

**[KLSBook]**
Roelof Koekoek , Peter A. Lesky , René F. Swarttouw, Hypergeometric Orthogonal Polynomials and Their q-Analogues, Springer 2010.

**[KS98Book]**
Leonid I. Korogodski, Yan S. Soibelman, Algebras of Functions on Quantum Groups: Part I, American Mathematical Soc., 1998

**[KS97]**
K. Krebs and S. Sandow, Matrix Product Eigenstates for One-Dimensional Stochastic Models and Quantum Spin Chains, J. Phys. A 30 (1997) 3165 and arXiv:cond-mat/9610029.

**[Kua18]**
J. Kuan, *International Mathematics Research Notices*, Volume 2018, Issue 17, September 2018, Pages 5378–5416.

**[Kua22]**
J. Kuan, *International Mathematics Research Notices*, Volume 2022, Issue 13, July 2022, Pages 9633–9662.

**[KZ24]**
Kuan, J. and Zhou, Z. (2024), Dualities and Asymptotics of Dynamic Stochastic Higher Spin Vertex Models, using the Drinfeld Twister,*arXiv preprint arXiv:2305.17602*.
https://arxiv.org/abs/2305.17602

**[Kua18]**
J. Kuan, An algebraic construction of duality functions for the stochastic $U_q(A_n^{(1)})$ vertex model and its degenerations, Commun. Math. Phys. (2018) 359: 121

**[KS92]**
P. P. Kulish and E. K. Sklyanin, *Algebraic structures related to reflection equations*, J. Phys. A: Math. Gen. **25** (1992), 5963–5975.

**[KSS93]**
P. P. Kulish, R. Sasaki and C. Schwiebert, *Constant solutions of reflection equations and quantum groups*, J. Math. Phys. **34** (1993), 286–304.

**[KMMO16]**
A. Kuniba, V.V. Mangazeev, S. Maruyama, M. Okado, Stochastic R matrix for $U_q(A_n^{(1)})$, Nuclear Physics B, Volume 913, 2016, Pages 248-277.

**[Li93]**
You-Quan Li, Yang Baxterization, J. Math. Phys. 34 (1993) 757.

**[Mac89]**
A.J. Macfarlane, On q-analogues of the quantum harmonic oscillator and the quantum group $SU(2)_q$. J. Phys. A, Volume 22, Number 21, 1989.

**[MGP68]**
MacDonald CT, Gibbs JH, Pipkin AC. Kinetics of biopolymerization on nucleic acid templates. Biopolymers. 1968;6(1):1–5.

**[Man14]**
Vladimir V. Mangazeev, On the Yang–Baxter equation for the six-vertex model, Nuclear Physics B, Volume 882, May 2014, Pages 70-96

**[RP81]**
L.C.G. Rogers and J.W. Pitman, Markov Functions, The Annals of Probability, Vol. 9, No. 4 (Aug., 1981), pp. 573-582 (10 pages)

**[Pau35]**
Pauling, L. (1935). "The Structure and Entropy of Ice and of Other Crystals with Some Randomness of Atomic Arrangement". Journal of the American Chemical Society. 57 (12): 2680–2684.

**[Pov13]**
A.M. Povolotsky, On the integrability of zero-range chipping models with factorized steady states J. Phys. A, Math. Theor., 46 (2013), p. 465205

**[San94]**
S. Sandow, Partially asymmetric exclusion process with open boundaries, Phys Review E, Volume 50, Number 4, 1994.

**[Sch97]**
G. Schütz, Duality relations for asymmetric exclusion processes, J Stat Phys 86, 1265–1287 (1997).

**[Sch97a]**
G. Schütz, Exact solution of the master equation for the asymmetric exclusion process. J. Stat Physics, **88** 427–445 (1997).

**[Spi70]**
F. Spitzer, Interaction of Markov Processes, Advances in Mathematics, Volume 5, Issue 2, October 1970, Pages 246–290.

**[Skl88]**
E. K. Sklyanin, *Boundary conditions for integrable quantum systems*, J. Phys. A: Math. Gen. **21** (1988), 2375–2389.

**[TW93]**
Tracy, Craig A. Widom, Harold, "Level-spacing distributions and the Airy kernel", Physics Letters B. 305 (1—2): 115-–118 (1993).

**[TW94]**
Tracy, C. A.; Widom, H. (1994), "Level-spacing distributions and the Airy kernel", Communications in Mathematical Physics, 159 (1): 151-–174 (1994).

**[TW08]**
C. Tracy and H. Widom, "Integral Formulas for the Asymmetric Simple Exclusion Process," Comm. Math. Phys., volume 279, pages 815–844 (2008).

**[TW09]**
Tracy, C. A.; Widom, H., "Asymptotics in ASEP with step initial condition", Communications in Mathematical Physics, 290 (1): 129–154, 2009.

**[ZGB91]**
R.B. Zhang, M.D. Gould and A.J. Bracken, From representations of the braid group to solutions of the Yang–Baxter equation, Nucl. Phys. B 354 (1991) 625.

**[ZG93]**
Yao-Zhong Zhang, Mark D. Gould, Quantum Affine Algebras and Universal R-Matrix with Spectral Parameter. arXiv:hep-th/9307007